\documentclass{IEEEtran}
\usepackage{cite}
\usepackage{amsmath,amssymb,amsfonts}
\usepackage{mathrsfs}
\usepackage{graphicx}
\usepackage{textcomp}
\usepackage{tikz}
\usetikzlibrary{arrows.meta,calc}
\usepackage[hidelinks]{hyperref}
\hypersetup{
  pdftitle={Adaptive Output Regulation of Linear Systems with Unknown Harmonic Exosystems and Additive Noise},
  pdfauthor={Bingsen Li and George Weiss}
}

\newtheorem{assumption}{Assumption}
\newtheorem{proposition}{Proposition}
\newtheorem{lemma}{Lemma}
\newtheorem{theorem}{Theorem}
\newtheorem{problem}{Problem}
\newtheorem{remark}{Remark}
\newtheorem{algorithm}{Algorithm}
\newcommand{\dd}{\mathrm d}
\newcommand{\ii}{\mathrm i}
\newcommand{\E}{\mathbb E}
\newcommand{\Prob}{\mathbb P}
\DeclareMathOperator{\col}{col}

\DeclareMathOperator{\rank}{rank}
\DeclareMathOperator{\Rea}{Re}

\begin{document}
\title{Adaptive Output Regulation of Linear Systems with Unknown Harmonic Exosystems and Additive Noise}
\author{Bingsen Li and George Weiss%
\thanks{The authors are with the School of Electrical
and Computer Engineering, Tel Aviv University,
Tel Aviv 6997801, Israel
(bingsenli9@gmail.com; gweiss@tauex.tau.ac.il).}}
\maketitle

\begin{abstract}
This paper addresses error-feedback output regulation for known
multivariable linear plants subject to an unknown exosystem and
additive Gaussian process and measurement noise. The exosystem generates constant
and sinusoidal signals, but its dimension, frequencies, and coupling
matrices are unavailable. An output-injection transformation
represents the external action through known channels. A causal
observer then produces a residual independent of the applied input,
whose filtered output separates into a harmonic signal, stationary
colored Gaussian noise, and a decaying transient. Using this observation model,
a Fourier procedure consistently identifies the number
and values of the active positive frequencies in the presence of  noise. The identified frequencies
determine a regulator through interpolation with stable residual
filters.
A contractive filter realization, nonexpansive state resizing, and
bounded gain updates ensure closed-loop moment bounds uniformly
over the selected filter orders. If the gain bound contains the ideal
interpolation gain, the adaptive plant response and tracking error
converge in mean square to those of the corresponding known-frequency
regulator. The limiting mean-square tracking error is a 
linear function of the process- and measurement-noise covariances.
In the absence of both noises, the tracking error converges to zero. 
\end{abstract}

\begin{IEEEkeywords}
Adaptive output regulation, unknown exosystem, tracking-error
feedback, frequency identification, process noise, measurement noise.
\end{IEEEkeywords}

\section{Introduction}

Accurate reference tracking and rejection of persistent disturbances
are important in motion control, vibration suppression, and power
conversion~\cite{LandauConstantinescuRey2005,LandauEtAl2011,
TeodorescuEtAl2006}. Constant and sinusoidal components provide useful
models for references and disturbances at a fixed operating condition,
and periodic tracking has also motivated repetitive
control~\cite{HaraEtAl1988}. In practice, the frequencies, number of
active components, and coupling channels may be unavailable.
Random disturbances may act on the plant, while only noisy
tracking-error measurements are available for feedback.

Output regulation concerns reference tracking and disturbance
rejection when the external signals are generated by an autonomous
exosystem. The aim is to stabilize the closed loop while driving
the tracking error to zero. For linear multivariable systems,
the regulator equations describe the required steady-state
response~\cite{Francis1977}, while the internal model principle
characterizes the controller dynamics needed for robust
regulation~\cite{FrancisWonham1976,Davison1976}.
These results form the basis of standard linear regulator
designs~\cite{Wonham1985,Trentelman2001} and extend to nonlinear
output regulation~\cite{IsidoriByrnes1990,
ByrnesPriscoliIsidori1997,Huang2004}.
Related regulation results have also been established for regular
infinite-dimensional systems~\cite{StateFeedbackRegularSystems,
RobustRegularInfiniteExosystems}.

When the external dynamics are unknown, adaptive sinusoidal
compensation and adaptive internal models can be used to address
uncertain frequencies~\cite{BodsonDouglas1997,
MarinoSantosuossoTomei2003,MarinoTomei2003,
SerraniIsidoriMarconi2001,Ding2003}.
Unknown and undermodeled exosystem orders have also been studied
for minimum-phase systems~\cite{MarinoTomei2007,MarinoTomei2011}.
Related disturbance-rejection methods address slowly varying
unknown sinusoidal parameters in MIMO linear time-varying
plants~\cite{SinusoidalMIMOLTV}.

The identification-based design in~\cite{BinMarconiTeel2019}
combines a continuous-time internal model with a discrete-time
identifier. Adaptation is formulated as a prediction problem,
separating the design of the predictor from stabilization.
A recursive least-squares scheme updates the internal-model
parameters from closed-loop data. Between parameter updates,
the temporary steady-state response retains the external frequencies
even when the estimated parameters are incorrect. This permits
recovery of the external dynamics under sufficient internal-model
dimension and persistent excitation.
The framework accommodates non-square, non-minimum-phase plants
and establishes robustness to additional unmodeled disturbances.
Related nonlinear formulations use continuous-time
identifiers~\cite{BernardBinMarconi2020} or relate approximate
regulation to the prediction error~\cite{BinBernardMarconi2021}.

Adaptive observers provide related state and parameter estimation
mechanisms~\cite{Kreisselmeier1977,Zhang2002}, including
observer-based error-feedback regulation for non-smooth,
non-periodic signals~\cite{NiuChenScarciotti2024}.
Stochastic output regulation has been investigated in
full-information and output-feedback formulations, including
systems with unstable linear
exosystems~\cite{Scarciotti2018,MelloneScarciotti2022,
MaYuLanChen2025}. Performance evaluation under regulation
constraints is also relevant to exosystem-observer
$H^2$ control~\cite{Hashikura2019}.
Here the exosystem is deterministic and unknown, while the plant
and measurements contain additive Gaussian noise.
Cancellation of the deterministic harmonics can therefore leave
a nonzero tracking variance, which must be distinguished from
the mean tracking error.

Frequency estimation is well developed through likelihood and
Fourier methods~\cite{RifeBoorstyn1974,QuinnFernandes1991,
QuinnHannan2001}, subspace and multichannel
techniques~\cite{Schmidt1986,RoyKailath1989,StoicaMoses2005,
GrantQuinn2020}, and statistical order
selection~\cite{WaxKailath1985,NadlerKontorovich2011}.
For error-feedback regulation, however, the external signal is
not directly measured. It must be extracted from a tracking error
that also depends on the plant response and applied input.
Moreover, stable filtering of white measurement noise produces
temporally correlated noise~\cite{SarkkaSolin2019}.
With an unknown frequency count, identification and control must
also remain well defined before the selected count becomes correct.
The problem considered here is therefore to identify the active
frequency count and values consistently from noisy error feedback,
while maintaining closed-loop moment bounds during adaptation.

This paper addresses this problem for known finite-dimensional
multivariable linear plants with possible direct feedthrough,
an unknown harmonic exosystem with fixed frequencies,
and additive Gaussian process and measurement noise.
The plant is assumed to be stabilizable and detectable and to
satisfy the standard regulation rank condition at the external
frequencies; the plant need not be stable or minimum phase.
The main contributions are as follows.
\begin{enumerate}
\item An output-injection transformation constructs an equivalent
system in which the unknown external action enters through known
channels and preserves the relation between the original and
transformed regulator equations. A causal observer uses the measured
error and the internally available input to generate a residual
whose dynamics are independent of the applied input.
Its filtered output separates into harmonics, stationary colored
Gaussian noise, and a decaying transient. Consequently, the
identification model remains valid throughout controller adaptation.

\item A growing-block Fourier procedure consistently identifies both
the number and values of the active positive frequencies in the
presence of process and measurement noise. Almost surely, the
selected count is eventually correct and the frequency error
multiplied by the block length tends to zero. For a prescribed band, the estimator requires no upper bound on
the frequency count, known lower bounds on amplitudes or separation,
or noise-covariance information; no small-noise condition is imposed.
It recovers the frequencies of the equivalent external signal
needed for regulation.

\item The identified frequencies determine an interpolation gain
for a stable residual-filter cascade, whose order is updated from
the estimated count. A uniformly contractive filter realization, nonexpansive resizing, and
bounded gain updates keep the closed-loop second and fourth
moments uniformly bounded, even while the frequency count
is still being learned.

\item Under the gain-budget condition, the adaptive tracking error
approaches that of the known-frequency regulator driven by the
same noise in mean square. Its mean tends to zero, and a Lyapunov
equation gives the limiting mean-square error from both noise
covariances. Exact asymptotic regulation is recovered when
both noises vanish.
\end{enumerate}

Sections~\ref{sec:problem}--\ref{sec:equivalent} formulate the problem
and construct the equivalent system and residual.
Sections~\ref{sec:identification}--\ref{sec:controller} develop
frequency identification and known-frequency regulation.
Section~\ref{sec:closed-loop} establishes adaptive closed-loop
performance, and Section~\ref{sec:numerical-simulations} presents
numerical results.

\emph{Notation:} The Euclidean norm 
and its induced matrix norm are denoted by $\|\cdot\|$, and
$\|\cdot\|_F$ is the Frobenius norm. 
For a finite set $\mathcal S$, $|\mathcal S|$ denotes its cardinality.
The spectrum of a square matrix $M$ is $\sigma(M)$; all ranks
at complex arguments are taken over $\mathbb C$. Column stacking
is denoted by $\col$, and $\otimes$ denotes the Kronecker product.
For a state trajectory $x$, $x(t^-):=\lim_{s\uparrow t}x(s)$ and
$x(t^+):=\lim_{s\downarrow t}x(s)$ denote its left and right limits.
Almost-sure statements are abbreviated as a.s. The notation
$X_k=o(a_k)$ a.s. means $\|X_k\|/a_k\to0$ as $k\to\infty$
on an event of probability one. 

\section{Problem Formulation}
\label{sec:problem}
Consider the finite-dimensional linear system
\begin{equation}\label{eq:original-system-feedthrough}
\begin{aligned}
 \dd z&=(Az+Bu+Pw)\,\dd t+V_p\,\dd W_p,\\
 \dot w&=Sw,\\
 e&=Cz+Du+Qw+e_{\mathrm{m}},
\end{aligned}
\end{equation}
where $z\in\mathbb R^n$, $u\in\mathbb R^m$,
$e\in\mathbb R^p$, and $w\in\mathbb R^{n_w}$. The initial states $z(0)$ and $w(0)$ are finite and deterministic.
The process-noise matrix $V_p\in\mathbb R^{n\times n_p}$ is
constant, and $W_p$ is a standard $n_p$-dimensional Wiener process.
Its intensity $\Sigma_p=V_pV_p^\top\succeq0$ need not be known.
The plant matrices $A,B,C,D$ are known. The exosystem matrix $S$,
the coupling matrices $P,Q$, and the state $w$ are \textit{unavailable}.
Only the tracking-error measurement $e$ and the internally available
input $u$ can be used by the controller. The term $e_{\mathrm{m}}$ represents
measurement noise; the true tracking error, before this noise is
added, is
\begin{equation}\label{eq:true-error}
 e_0=Cz+Du+Qw.
\end{equation}

\begin{assumption}\label{ass:plant-exosystem}
The pair $(A,C)$ is detectable and the pair $(A,B)$ is stabilizable. The real matrix $S$ is diagonalizable
over $\mathbb C$ and satisfies $\sigma(S)\subset\ii\mathbb R$.
\end{assumption}

Under this assumption, $w(t)$ is bounded and consists of constant
and sinusoidal components. No value of $n_w$ is needed in the
implementation. To give a precise meaning to Gaussian white
measurement noise, replace the last equation of
\eqref{eq:original-system-feedthrough} by the measurement process
\begin{equation}\label{eq:measurement-process}
 \dd E(t)=e_0(t)\,\dd t+V_m\,\dd W_m(t),\qquad E(0)=0,
\end{equation}
where $W_m$ is a standard $p$-dimensional Wiener process, independent
of $W_p$, and $V_m\in\mathbb R^{p\times p}$ is constant.
All stochastic differential equations are interpreted in the It\^o sense. The noise intensity
$\Sigma_m=V_mV_m^\top\succeq0$ need not be known.
The formal notation is $e=\dot E$ and $e_{\mathrm{m}}=V_m\dot W_m$;
all implemented equations use $\dd E$, not pointwise white-noise
values.

For $s\in\mathbb C$, define the control-channel Rosenbrock matrix
\begin{equation}\label{eq:Ru}
 \mathcal R_u(s)=\begin{bmatrix}sI-A&-B\\C&D\end{bmatrix}.
\end{equation}
The standard regulation feasibility condition
\cite[Ch.~9]{Trentelman2001} is imposed as follows.

\begin{assumption}\label{ass:regulation}
The Rosenbrock matrix \eqref{eq:Ru} is full rank at $\sigma(S)$, i.e., 
\begin{equation}\label{eq:Ru-rank}
 \rank\mathcal R_u(\lambda)=n+p,
 \qquad\lambda\in\sigma(S).
\end{equation}
\end{assumption}

\begin{problem}\label{prob:adaptive-regulation}
Using the noisy observation
\eqref{eq:measurement-process} and the known plant matrices,
achieve the following objectives:
\begin{enumerate}
 \item\label{obj:identification}
 Consistently identify the number and values of the active
 frequencies required for regulation.

 \item\label{obj:regulation}
 Design an error-feedback controller such that the closed loop
 has a unique global strong solution, the plant and controller
 states and the control input remain bounded in mean square,
 and the tracking error satisfies
 \begin{equation}\label{eq:regulation-objective}
 \begin{aligned}
  \E[e_0(t)]&\rightarrow0,\\
  \limsup_{t\to\infty}\E\|e_0(t)\|^2
  &\leq c\bigl[\operatorname{tr}(\Sigma_m)+\operatorname{tr}(\Sigma_p)\bigr],
 \end{aligned}
 \end{equation}
 where $c\geq0$ is independent of $\Sigma_m$ and $\Sigma_p$.
\end{enumerate}
\end{problem}
\section{Equivalent System and Residual Construction}
\label{sec:equivalent}
The aim of this section is to obtain a measured signal with a harmonic
component and no dependence on the applied input. We first represent
the unknown external action through known injection matrices. We
then construct a residual from $(E,u)$ and characterize the transient
and noise that accompany its harmonic part.

\subsection{Representing the external action through known channels}
By detectability, choose $L_o\in\mathbb R^{n\times p}$ such that
\begin{equation}\label{eq:Ao}
 A_o:=A-L_oC
\end{equation}
is Hurwitz. This gain is fixed throughout the construction.
The spectra of $A_o$ and $S$ are disjoint, so the standard Sylvester
equation result~\cite{Trentelman2001} gives a unique real
$\Pi\in\mathbb R^{n\times n_w}$ satisfying
\begin{equation}\label{eq:sylvester-equivalent}
 \Pi S=A_o\Pi+P-L_oQ.
\end{equation}
Define the equivalent plant state $\bar z$, the external-signal
matrix $\Gamma\in\mathbb R^{p\times n_w}$, and the signal
$d\in\mathbb R^p$ by
\begin{equation}\label{eq:equivalent-definitions}
 \bar z=z-\Pi w,\qquad \Gamma=C\Pi+Q,\qquad d=\Gamma w.
\end{equation}
Neither $\Pi$ nor $\Gamma$ is required by the implementation.
They serve to identify the signal represented by the residual.

Indeed, \eqref{eq:sylvester-equivalent} gives
$A\Pi+P-\Pi S=L_o\Gamma$. Substituting
\eqref{eq:equivalent-definitions} into the state equation
\eqref{eq:original-system-feedthrough} and the measurement equation \eqref{eq:measurement-process}, and using
\eqref{eq:sylvester-equivalent} together with
$d=\Gamma w=(C\Pi+Q)w$, yields the \textit{equivalent system}:
\begin{equation}\label{eq:equivalent-plant}
\begin{aligned}
 \dd\bar z&=(A\bar z+Bu+L_od)\,\dd t+V_p\,\dd W_p,\\
 \dot w&=Sw,\qquad d=\Gamma w,\\
 \dd E&=e_0 \,\dd t+V_m\,\dd W_m, \qquad e_0=C\bar z+Du+d.
\end{aligned} 
\end{equation}
The transformation is invertible because $z=\bar z+\Pi w$.
Thus the unknown channels $P,Q$ have been replaced by the known
injection matrices $L_o, I_p$. In particular, the Rosenbrock matrix
from $d$ to the true error is
\begin{equation}\label{eq:Rd}
 \mathcal R_d(s)=\begin{bmatrix}sI-A&-L_o\\C&I_p\end{bmatrix},
 \qquad\det\mathcal R_d(s)=\det(sI-A_o).
\end{equation}
Clearly, it is nonsingular on the imaginary axis.

Under Assumption~\ref{ass:regulation}, the equivalent system admits
 $X\in\mathbb R^{n\times n_w}$ and
$U\in\mathbb R^{m\times n_w}$ satisfying the regulator
equations (see, for instance,
\cite{Francis1977,Trentelman2001})
\begin{equation}\label{eq:equivalent-regulator-equations}
 XS=AX+BU+L_o\Gamma,\qquad CX+DU+\Gamma=0.
\end{equation}
Using \eqref{eq:sylvester-equivalent} and $\Gamma=C\Pi+Q$,
we find that $X_{\mathrm{ori}}=\Pi+X$ satisfies the regulator equations
of the original system:
\begin{equation}\label{eq:original-regulator-equations}
 X_{\mathrm{ori}}S=AX_{\mathrm{ori}}+BU+P,\quad CX_{\mathrm{ori}}+DU+Q=0.
\end{equation}
Thus, in the absence of both noises, $\bar z=Xw$ and $z=X_{\mathrm{ori}}w$ are corresponding
steady-state responses with the same input $u=Uw$ and zero
tracking error. These equations are used only for analysis;
the controller below uses the identified frequencies and
does not require $X$ or $U$.

\subsection{An implementable residual and its filtered output}
To estimate $\bar z$,
introduce an auxiliary observer state $\xi\in\mathbb R^n$:
\begin{equation}\label{eq:implementable-observer}
 \dd\xi=\bigl[A_o\xi+(B-L_oD)u\bigr]\dd t+L_o\,\dd E
\end{equation}
from a prescribed finite deterministic $\xi(0)$. Subtract its
predicted output from the measurement to form the integrated
residual $\mathsf R\in\mathbb R^p$:
\begin{equation}\label{eq:integrated-residual}
 \dd\mathsf R=\dd E-(C\xi+Du)\,\dd t,\qquad\mathsf R(0)=0.
\end{equation}
Both equations use only $(E,u)$ and known matrices. To determine
what $\mathsf R$ contains, define 
\begin{equation}\label{eq:chi}
 \chi=\bar z-\xi=z-\xi-\Pi w.
\end{equation}
Subtracting \eqref{eq:implementable-observer} from
\eqref{eq:equivalent-plant} gives
\begin{equation}\label{eq:residual-system}
\begin{aligned}
 \dd\chi&=A_o\chi\,\dd t-L_oV_m\,\dd W_m+V_p\,\dd W_p,\\
 \dd\mathsf R&=(d+C\chi)\,\dd t+V_m\,\dd W_m.
\end{aligned}
\end{equation}
The applied input has canceled. This identity remains valid when
$u$ is generated by feedback from the same measurement, and will
allow the estimator to operate throughout controller adaptation.

To obtain an continuous signal for Fourier analysis,
choose $\alpha>0$ and filter $\mathsf R$ as
\begin{equation}\label{eq:lowpass-implementation}
 \dd y_f=-\alpha y_f\,\dd t+\alpha\,\dd\mathsf R,
 \qquad y_f(0)\in\mathbb R^p,
\end{equation}
where $y_f(0)$ is finite and deterministic. The corresponding
scalar transfer function and its vector version are
\begin{equation}\label{eq:lowpass-transfer}
 h_f(s)=\frac{\alpha}{s+\alpha},\qquad H_f(s)=h_f(s)I_p.
\end{equation}
Since $h_f(0)=1$ and $h_f(\ii\omega)\ne0$ for every real
$\omega$, this filter preserves each nonzero harmonic component
present in $d$. Its effect on noise is examined next.

\subsection{Separating the harmonic signal, noise, and transient}
The steady-state response of \eqref{eq:lowpass-implementation} to
$d=\Gamma w$ is $\Gamma_fw$, where
\begin{equation}\label{eq:Gammaf}
 \Gamma_f=\alpha\Gamma(S+\alpha I)^{-1},
 \qquad \Gamma_fS=-\alpha\Gamma_f+\alpha\Gamma.
\end{equation}
Subtract this response and collect the remaining filter state with
the observer error in $\zeta=\col(\chi,y_f-\Gamma_fw)$. Define
\begin{equation}\label{eq:noise-realization}
\begin{gathered}
 \mathcal A=\begin{bmatrix}A_o&0\\\alpha C&-\alpha I_p\end{bmatrix},
 \quad\mathcal B_m=\begin{bmatrix}-L_o\\\alpha I_p\end{bmatrix},\\
 \mathcal C=\begin{bmatrix}0&I_p\end{bmatrix},\quad
 \mathcal B_p=\begin{bmatrix}I_n\\0_{p\times n}\end{bmatrix}.
\end{gathered}
\end{equation}
Then \eqref{eq:residual-system}--\eqref{eq:Gammaf} give
\begin{equation}\label{eq:zeta-dynamics}
\begin{aligned}
 \dd\zeta&=\mathcal A\zeta\,\dd t+\mathcal B_mV_m\,\dd W_m
                   +\mathcal B_pV_p\,\dd W_p,\\
 y_f&=\Gamma_fw+\mathcal C\zeta.
\end{aligned}
\end{equation}
Since $\mathcal A$ is Hurwitz, a stationary Gaussian solution
is~\cite[Sec.~6.5]{SarkkaSolin2019}
\begin{equation}\label{eq:stationary-state}
\begin{aligned}
 \zeta_s(t)={}&\int_{-\infty}^t
 e^{\mathcal A(t-s)}\mathcal B_mV_m\,\dd W_m(s)\\
 &+\int_{-\infty}^t
 e^{\mathcal A(t-s)}\mathcal B_pV_p\,\dd W_p(s).
\end{aligned}
\end{equation}
where $W_m$ and $W_p$ are extended to negative times using mutually
independent Wiener histories, independent also of their future increments.
The extension does not alter any measurement for $t\geq0$.
Because $\zeta$ and $\zeta_s$ have the same driving increments
after time zero,
\begin{equation}\label{eq:stationary-coupling}
 \zeta(t)-\zeta_s(t)=e^{\mathcal At}[\zeta(0)-\zeta_s(0)].
\end{equation}
Based on \eqref{eq:original-system-feedthrough}, \eqref{eq:zeta-dynamics} and \eqref{eq:stationary-coupling}, the measured filter output has the exact decomposition
\begin{equation}\label{eq:filtered-decomposition}
 y_f(t)=\Gamma_fe^{St}w(0)+\nu_f(t)+r_f(t),
 \qquad t\geq0,
\end{equation}
where $\nu_f=\mathcal C\zeta_s$ is the stationary noise and
$r_f(t)=\mathcal Ce^{\mathcal At}[\zeta(0)-\zeta_s(0)]$
is the initialization transient. In particular, for a deterministic
$\beta>0$ and an a.s. finite random $M_{\mathrm{tr}}$,
\begin{equation}\label{eq:transient-decay}
 \|r_f(t)\|\leq M_{\mathrm{tr}}e^{-\beta t},\qquad t\geq0.
\end{equation}

Since $\mathcal A$ is Hurwitz, the stationary covariance
$\Sigma_{\zeta}$ of $\zeta_s$ is the unique positive semidefinite
solution of the Lyapunov equation, see for instance,~\cite[Sec.~6.5]{SarkkaSolin2019},
\begin{equation}\label{eq:stationary-covariance}
\begin{aligned}
 &\mathcal A\Sigma_{\zeta}+\Sigma_{\zeta}\mathcal A^\top
 +\mathcal B_m\Sigma_m\mathcal B_m^\top+\mathcal B_p\Sigma_p\mathcal B_p^\top=0.
\end{aligned}
\end{equation}
For $\tau\geq0$, we know
\[
\begin{aligned}
 \zeta_s(t+\tau)={}&e^{\mathcal A\tau}\zeta_s(t)
 +\int_t^{t+\tau}
 e^{\mathcal A(t+\tau-s)}\mathcal B_mV_m\,\dd W_m(s)\\
 &+\int_t^{t+\tau}
 e^{\mathcal A(t+\tau-s)}\mathcal B_pV_p\,\dd W_p(s).
\end{aligned}
\]
Both stochastic integrals have zero mean and are independent of
$\zeta_s(t)$. Since $\zeta_s$ has zero mean and stationary
covariance $\Sigma_{\zeta}$, it follows that
\[
 \E[\zeta_s(t+\tau)\zeta_s(t)^\top]
 =e^{\mathcal A\tau}\Sigma_{\zeta}.
\]
Hence, using $\nu_f=\mathcal C\zeta_s$, we obtain
\begin{equation}\label{eq:noise-covariance}
\begin{aligned}
 &R_\nu(\tau)
 :=\E[\nu_f(t+\tau)\nu_f(t)^\top]=\mathcal Ce^{\mathcal A\tau}\Sigma_{\zeta}\mathcal C^\top,
 \quad \tau\geq0,\\
 &R_\nu(-\tau)=R_\nu(\tau)^\top.
\end{aligned}
\end{equation}
Exponential decay of $e^{\mathcal A\tau}$ implies
\begin{equation}\label{eq:cov-integrable}
 K_\nu:=\int_{\mathbb R}\|R_\nu(\tau)\|\,\dd\tau<\infty.
\end{equation}

Thus, \eqref{eq:filtered-decomposition} separates $y_f$ into
a harmonic signal, zero-mean stationary Gaussian noise with
integrable covariance, and an exponentially decaying transient.
These properties will be used in the identification analysis
of Section~\ref{sec:identification}. The sum $\nu_f+r_f$ need
not be stationary, and independence between $\nu_f$ and $r_f$
is not required.

\subsection{Frequencies available for identification}
The deterministic part of \eqref{eq:filtered-decomposition} is
harmonic because $S$ is diagonalizable with imaginary-axis spectrum.
To specify which frequencies are present, let $\mathcal P_\lambda$
be the spectral projector of $S$ at $\lambda$, with
$\mathcal P_0=0$ when $0\notin\sigma(S)$. Define the set of
active frequencies by
\begin{equation}\label{eq:active-frequencies}
 \Omega_d=\{\omega>0:\ii\omega\in\sigma(S),\ 
                    \Gamma\mathcal P_{\ii\omega}w(0)\ne0\}.
\end{equation}
Let $q=|\Omega_d|$. When $q>0$, write its distinct elements as
$0<\omega_1<\cdots<\omega_q$. Then
\begin{equation}\label{eq:harmonic-d}
\begin{aligned}
 d(t)&=d_0+\sum_{j=1}^q\bigl(c_je^{\ii\omega_jt}
                         +\overline{c_j}e^{-\ii\omega_jt}\bigr),\\
 d_0&=\Gamma\mathcal P_0w(0),\qquad
 c_j=\Gamma\mathcal P_{\ii\omega_j}w(0)\ne0.
\end{aligned}
\end{equation}
Here $d_0\in\mathbb R^p$ and $c_j\in\mathbb C^p$.
The case $q=0$ is allowed, and $d_0$ may vanish. Combining
\eqref{eq:Gammaf} with \eqref{eq:harmonic-d} gives the observation
model used in the next section:
\begin{equation}\label{eq:identification-model}
\begin{aligned}
 y_f(t)={}&d_0+\sum_{j=1}^q
 \bigl(c_{j,f}e^{\ii\omega_jt}
      +\overline{c_{j,f}}e^{-\ii\omega_jt}\bigr)\\
 &+\nu_f(t)+r_f(t),\qquad
 c_{j,f}=h_f(\ii\omega_j)c_j\ne0.
\end{aligned}
\end{equation}

\begin{remark}[Scope of frequency recovery]\label{rem:visibility}
The target is the active set $\Omega_d$, not necessarily the entire
spectrum of $S$. An exosystem mode absent from $w(0)$, or removed
by $\Gamma$, does not appear in $d$. The control-channel rank
condition \eqref{eq:Ru-rank} does not prevent this absence.
For regulation of the given signal, only the frequencies present
in $d$ need to be considered.  The constant component is removed for positive-frequency
identification and is handled separately by the controller.
\end{remark}

\section{Online Identification of the Equivalent Signal}
\label{sec:identification}
The residual construction has reduced identification to the harmonic
observation model \eqref{eq:identification-model}. The controller
will need the number $q$ of active positive frequencies and their
values $\omega_1,\ldots,\omega_q$. To estimate these quantities,
we accumulate Fourier coefficients over successive observation
blocks and select peaks above a decreasing threshold. A blockwise
mean removes the constant component before this search. 
\subsection{Search interval and observation blocks}
\begin{assumption}[Search interval]\label{ass:band}
Known constants $0<\omega_-<\omega_+$ satisfy
$\Omega_d\subset(\omega_-,\omega_+)$. The frequencies and
coefficients in \eqref{eq:harmonic-d} are fixed in time.
\end{assumption}
The estimator requires neither the number of frequencies nor
lower bounds on their amplitudes or separation. Consistency is
established for each fixed signal with distinct active frequencies.
The controller handles zero frequency separately.

Choose $t_{\mathrm{b}}>0$ and $\gamma_{\mathrm{b}}>0$, and set
$\Delta\omega=\omega_+-\omega_-$. For $k=0,1,\ldots$, define
\begin{equation}\label{eq:blocks}
\begin{gathered}
 t_k=t_{\mathrm{b}}2^k,\ \
 T_k=t_{k+1}-t_k=t_{\mathrm{b}}2^k,\ \
 I_k=[t_k,t_{k+1}].
\end{gathered}
\end{equation}
On block $I_k$, use the grid spacing $h_k$, exclusion parameter
$\rho_k$, and detection threshold $\gamma_k$:
\begin{equation}\label{eq:grid-scales}
 h_k=(\Delta\omega)2^{-2k},\
 \rho_k=(\Delta\omega)2^{-k/4},\
 \gamma_k=\gamma_{\mathrm{b}}2^{-k/4}.
\end{equation}
For each block $I_k$,  define the frequency search grid by
\begin{equation}\label{eq:frequency-grid}
 \mathcal G_k
 =\{\omega_-+\ell h_k:\ell=0,1,\ldots,2^{2k}\}.
\end{equation}
These choices do not require either noise intensity $\Sigma_m$ or $\Sigma_p$.
They satisfy, as $k\to\infty$,
\begin{equation}\label{eq:scale-properties}
 T_kh_k\to0,\qquad
 \rho_k,\gamma_k\to0,\qquad
 (T_k\rho_k)^{-1}=o(\gamma_k).
\end{equation}
\subsection{Blockwise frequency estimation}
On each block $I_k$, the estimator accumulates Fourier integrals
and selects frequency peaks at $t_{k+1}$. Introduce
$m_k\in\mathbb R^p$ and $Z_{k,\omega}\in\mathbb C^p$ with dynamics
\begin{equation}\label{eq:online-accumulators}
 \dot m_k=y_f(t),\qquad
 \dot Z_{k,\omega}=y_f(t)e^{-\ii\omega t},
 \quad \omega\in\mathcal G_k.
\end{equation}
Define the known scalar factor
\begin{equation}\label{eq:block-kernel}
 D_k(\nu)=\frac1{T_k}\int_{I_k}e^{\ii\nu t}\,\dd t,
 \qquad \nu\in\mathbb R.
\end{equation}
At the block endpoint, the mean and centered Fourier coefficients
are computed as
\begin{equation}\label{eq:mean-Fourier}
\begin{aligned}
 \bar y_k&=\frac{m_k(t_{k+1})}{T_k},\\
 F_k(\omega)
 &=\frac1{T_k}\int_{I_k}
       [y_f(t)-\bar y_k]e^{-\ii\omega t}\,\dd t\\
 &=\frac{Z_{k,\omega}(t_{k+1})}{T_k}
       -\bar y_kD_k(-\omega).
\end{aligned}
\end{equation}
The complete procedure is given
in Algorithm~\ref{alg:frequency-identification}.

\begin{algorithm}[Frequency estimation on block $I_k$]
\label{alg:frequency-identification}
\normalfont

\noindent\textbf{Step 1. Accumulate the data.}
At $t_k$, set $m_k=0$ and $Z_{k,\omega}=0$ for all
$\omega\in\mathcal G_k$. Integrate
\eqref{eq:online-accumulators} over $I_k$.

\par\smallskip
\noindent\textbf{Step 2. Compute the Fourier coefficients.}
At $t_{k+1}$, compute $\bar y_k$ and $F_k(\omega)$ using
\eqref{eq:mean-Fourier}. Initialize the remaining candidate set
$\mathcal G_k^{\mathrm{rem}}=\mathcal G_k$ and the selected set
$\widehat\Omega_k=\varnothing$.

\par\smallskip
\noindent\textbf{Step 3. Find the largest remaining peak.}
If $\mathcal G_k^{\mathrm{rem}}=\varnothing$, go to Step~5. Otherwise, select
\begin{equation}\label{eq:peak-search}
 \widehat\omega\in
 \mathop{\arg\max}_{\omega\in\mathcal G_k^{\mathrm{rem}}}\|F_k(\omega)\|,
\end{equation}
using a fixed rule to resolve ties.
If $\|F_k(\widehat\omega)\|\leq\gamma_k$, go to Step~5.

\par\smallskip
\noindent\textbf{Step 4. Record the peak and exclude its neighbors.}
Add $\widehat\omega$ to $\widehat\Omega_k$ and replace
\begin{equation}\label{eq:peak-deletion}
 \mathcal G_k^{\mathrm{rem}}=
 \{\omega\in\mathcal G_k^{\mathrm{rem}}:|\omega-\widehat\omega|>2\rho_k\}.
\end{equation}
Return to Step~3.

\par\smallskip
\noindent\textbf{Step 5. Output the estimates.}
Set $\widehat q_k=|\widehat\Omega_k|$ and sort the selected
frequencies as
$\widehat\omega_{1,k}<\cdots<
 \widehat\omega_{\widehat q_k,k}$.
Release these estimates at $t_{k+1}$.
If no peak is selected, return $\widehat q_k=0$ and an empty list.

\end{algorithm}

The norm $\|F_k(\omega)\|$ combines information from all output
components. The threshold $\gamma_k$ rejects small peaks, while
the exclusion radius $2\rho_k$ avoids selecting several nearby
grid points for the same peak. The procedure is repeated on each
block. Only the estimated frequencies and their number are supplied
to the controller; the block mean is used only for centering.
\subsection{Uniform control of the colored-noise coefficients}
The grid grows with the observation length, so a pointwise noise
estimate at a fixed frequency does not suffice. The next lemma
controls every candidate frequency simultaneously. It applies to
any stationary Gaussian noise with integrable covariance; the
construction in Section~\ref{sec:equivalent} satisfies these
conditions by \eqref{eq:stationary-state}--\eqref{eq:cov-integrable}.

\begin{lemma}[Noise bound over the search grid]
\label{lem:noise-Fourier}
Let $\nu_f(t)\in\mathbb R^p$ be a jointly measurable, zero-mean,
finite-variance stationary Gaussian process whose covariance
$R_\nu$ satisfies \eqref{eq:cov-integrable}.
For the blocks, threshold, and grid in
\eqref{eq:blocks}--\eqref{eq:frequency-grid}, define
\begin{equation}\label{eq:noise-Fourier}
 N_k(\omega)=\frac1{T_k}\int_{I_k}
 \nu_f(t)e^{-\ii\omega t}\,\dd t.
\end{equation}
Then, as $k\to\infty$,
\begin{equation}\label{eq:uniform-noise}
 \max_{\omega\in\mathcal G_k\cup\{0\}}\|N_k(\omega)\|
 =o(\gamma_k)\quad\text{a.s.}
\end{equation}
\end{lemma}
\begin{IEEEproof}
Since $\nu_f$ is zero-mean and stationary, Fubini's theorem gives
\begin{align}
 \E|N_{k,a}(\omega)|^2
 &=\frac{1}{T_k^2}\int_{I_k}\int_{I_k}
 \E[\nu_{f,a}(t)\nu_{f,a}(s)]
 e^{-\ii\omega(t-s)}\,\dd t\,\dd s \notag\\
 &\leq\frac{1}{T_k^2}\int_{I_k}\int_{I_k}
 \|R_\nu(t-s)\|\,\dd t\,\dd s \notag\\
 &\leq\frac{1}{T_k^2}\int_{I_k}
 \left(\int_{\mathbb R}\|R_\nu(\tau)\|\,\dd\tau\right)
 \dd s =\frac{K_\nu}{T_k}.
 \label{eq:noise-variance}
\end{align}
Here the first inequality bounds each scalar covariance by
the matrix covariance norm, and the second follows by setting
$\tau=t-s$ and extending the inner integral to $\mathbb R$.

If $K_\nu=0$, \eqref{eq:noise-variance} gives $N_k(\omega)=0$
a.s. at all grid points and at zero, so \eqref{eq:uniform-noise}
is immediate. Hence assume $K_\nu>0$.

The real and imaginary parts of each component of $N_k(\omega)$
are Gaussian with mean zero and variance at most $K_\nu/T_k$.
Applying the Gaussian tail bound to each of these $2p$ real-valued
variables gives
\begin{equation}\label{eq:Gaussian-tail}
 \Prob\{\|N_k(\omega)\|>\varepsilon\}
 \leq4p\exp\!\left(-\frac{T_k\varepsilon^2}{4pK_\nu}\right).
\end{equation}
Take $\varepsilon_k=\gamma_{\mathrm{b}}2^{-3k/8}$. Applying \eqref{eq:Gaussian-tail} at each frequency and summing
over the $2^{2k}+2$ points in $\mathcal G_k\cup\{0\}$ gives
\begin{equation}\label{eq:BC-bound}
\begin{aligned}
 &\Prob\!\left\{
 \max_{\omega\in\mathcal G_k\cup\{0\}}
 \|N_k(\omega)\|>\varepsilon_k
 \right\}
\leq
 \sum_{\omega\in\mathcal G_k\cup\{0\}}
 \Prob\{\|N_k(\omega)\|>\varepsilon_k\}\\
 &\quad\leq
 4p(2^{2k}+2)
 \exp\!\left(-\frac{T_k\varepsilon_k^2}{4pK_\nu}\right)=
 C_1(2^{2k}+2)e^{-C_2 2^{k/4}},
\end{aligned}
\end{equation}
where $C_1=4p$ and $C_2=t_{\mathrm{b}}\gamma_{\mathrm{b}}^2/(4pK_\nu)>0$.
The last equality uses $T_k=t_{\mathrm{b}}2^k$. The right-hand side is summable in $k$. The first Borel--Cantelli
lemma implies that the maximum is at most $\varepsilon_k$ for all
sufficiently large $k$, a.s. Since
$\varepsilon_k/\gamma_k=2^{-k/8}\to0$, this proves
\eqref{eq:uniform-noise}. 
\end{IEEEproof}

\subsection{Consistency of the order and frequencies}
We next combine the noise bound with estimates of the transient
and finite-window leakage to prove that the algorithm identifies
the correct number and values of the active frequencies in \eqref{eq:identification-model}.
The following result states precisely which frequencies are
recovered and when the number of returned estimates becomes correct.

\begin{theorem}[Identification of the frequencies]
\label{thm:identification}
Suppose Assumptions~\ref{ass:plant-exosystem} and \ref{ass:band}
hold. Apply Algorithm~\ref{alg:frequency-identification} to
the observation \eqref{eq:identification-model}, using
\eqref{eq:grid-scales}.
For the active frequencies in \eqref{eq:harmonic-d}, there is,
on an event of probability one, a finite random integer $k_0$ such that
\begin{equation}\label{eq:order-consistency}
 \widehat q_k=q\qquad\text{for every }k\geq k_0,
\end{equation}
and
\begin{equation}\label{eq:frequency-consistency}
 T_k|\widehat\omega_{j,k}-\omega_j|\longrightarrow0
 \quad\text{as }k\to\infty,\qquad j=1,\ldots,q.
\end{equation}
\end{theorem}
\begin{IEEEproof}
 By \eqref{eq:transient-decay} and
Lemma~\ref{lem:noise-Fourier}, there is an event of probability
one on which
$\int_0^\infty\|r_f(t)\|\,\dd t<\infty$
and the noise estimate \eqref{eq:uniform-noise} holds. Fix a sample path in this event.
All limits below are taken along that path. Direct integration of \eqref{eq:block-kernel} gives
\begin{equation}\label{eq:kernel-bound}
\begin{aligned}
 D_k(\nu)
 &=e^{\ii\nu(t_k+T_k/2)}
   \frac{\sin(\nu T_k/2)}{\nu T_k/2},\\
 |D_k(\nu)|
 &\leq\min\{1,2/(T_k|\nu|)\},\qquad \nu\neq0.
\end{aligned}
\end{equation}
The bound follows from $|e^{\ii x}|=1$ and
$|\sin x|\leq\min\{|x|,1\}$.
The normalized transient coefficient $R_k(\omega)$ satisfies
\begin{equation}\label{eq:transient-Fourier}
\begin{aligned}
 R_k(\omega)&=\frac1{T_k}\int_{I_k}r_f(t)e^{-\ii\omega t}\,\dd t,\\
 \sup_{\omega\in\mathbb R}\|R_k(\omega)\|
 &\leq\frac{\|r_f\|_{L^1(0,\infty)}}{T_k}=o(\gamma_k).
\end{aligned}
\end{equation}

\emph{Centering:} Averaging \eqref{eq:identification-model} over $I_k$ and using
the definitions in \eqref{eq:mean-Fourier},
\eqref{eq:block-kernel}, \eqref{eq:noise-Fourier}, and
\eqref{eq:transient-Fourier} gives
\begin{equation}\label{eq:mean-decomposition}
\begin{split}
 \bar y_k-d_0
 ={}&\sum_{j=1}^q\bigl[c_{j,f}D_k(\omega_j)
               +\overline{c_{j,f}}D_k(-\omega_j)\bigr]\\
 &+N_k(0)+R_k(0).
\end{split}
\end{equation}
By \eqref{eq:kernel-bound},
the sum in \eqref{eq:mean-decomposition} is $O(T_k^{-1})$.
Lemma~\ref{lem:noise-Fourier} and
\eqref{eq:transient-Fourier} yield the estimate
\begin{equation}\label{eq:dc-consistency}
 \bar y_k-d_0=o(\gamma_k).
\end{equation}

\emph{Fourier decomposition:}
Substitute the signal model \eqref{eq:identification-model}
into the centered Fourier coefficient \eqref{eq:mean-Fourier}.
Using the definitions of $D_k$, $N_k$, and $R_k$ in
\eqref{eq:block-kernel}, \eqref{eq:noise-Fourier}, and
\eqref{eq:transient-Fourier}, respectively, gives, for
$\omega\in\mathcal G_k$,
\begin{equation}\label{eq:F-decomposition}
 F_k(\omega)=\sum_{j=1}^q c_{j,f}D_k(\omega_j-\omega)
                         +\mathcal E_k(\omega),
\end{equation}
where 
\begin{equation}\label{eq:E-decomposition}
\begin{split}
 \mathcal E_k(\omega)={}&\sum_{j=1}^q\overline{c_{j,f}}
                    D_k(-\omega_j-\omega)+N_k(\omega)\\
 &+R_k(\omega)+(d_0-\bar y_k)D_k(-\omega).
\end{split}
\end{equation}
For $\omega\in\mathcal G_k$, from 
\eqref{eq:frequency-grid}, we know $\omega\geq\omega_-$.
Since $\omega_j>\omega_-$, the kernel bound
\eqref{eq:kernel-bound} yields
\[
 |D_k(-\omega_j-\omega)|
 \leq\frac{2}{T_k(\omega_j+\omega)}
 \leq\frac{1}{T_k\omega_-}.
\]
Thus the negative-frequency sum in \eqref{eq:E-decomposition}
is uniformly $O(T_k^{-1})$.
By \eqref{eq:uniform-noise} and \eqref{eq:transient-Fourier},
$N_k(\omega)=o(\gamma_k)$ and $R_k(\omega)=o(\gamma_k)$
uniformly for $\omega\in\mathcal G_k$.
By \eqref{eq:dc-consistency} and \eqref{eq:kernel-bound},
$(d_0-\bar y_k)D_k(-\omega)=o(\gamma_k)$
uniformly for $\omega\in\mathcal G_k$.
Since \eqref{eq:blocks} and \eqref{eq:grid-scales} give
$T_k^{-1}=o(\gamma_k)$, we obtain
\begin{equation}\label{eq:uniform-remainder}
 \max_{\omega\in\mathcal G_k}\|\mathcal E_k(\omega)\|=o(\gamma_k).
\end{equation}

\emph{Localization and order selection:}
First suppose $q>0$.
For each $j\in\{1,\ldots,q\}$, let $g_{j,k}\in\mathcal G_k$
be a nearest grid point to $\omega_j$.
By \eqref{eq:frequency-grid},
$|g_{j,k}-\omega_j|\leq h_k$.
Since $T_k|\omega_j-g_{j,k}|\leq T_kh_k\to0$
by \eqref{eq:scale-properties} and $t_k=T_k$
by \eqref{eq:blocks}, \eqref{eq:kernel-bound} gives
$D_k(\omega_j-g_{j,k})\to1$.  For $\ell\neq j$, we have
$|\omega_\ell-g_{j,k}|\to|\omega_\ell-\omega_j|>0$.
Hence \eqref{eq:kernel-bound} gives
$D_k(\omega_\ell-g_{j,k})=O(T_k^{-1})$.
Also, \eqref{eq:uniform-remainder} and $\gamma_k\to0$
imply $\mathcal E_k(g_{j,k})\to0$.
Substituting these estimates into \eqref{eq:F-decomposition}
yields
\[
 F_k(g_{j,k})
 =c_{j,f}D_k(\omega_j-g_{j,k})+o(1).
\]
Together with the preceding limit for $D_k$, this proves
\begin{equation}\label{eq:nearest-peak}
 D_k(\omega_j-g_{j,k})\to1,
 \qquad F_k(g_{j,k})\to c_{j,f}.
\end{equation}

Using the radius $\rho_k$ from \eqref{eq:grid-scales},
define
\begin{equation}\label{eq:peak-neighborhoods}
 \mathcal U_{j,k}
 =\{\omega\in\mathcal G_k:|\omega-\omega_j|<\rho_k\}.
\end{equation}
Since $h_k/\rho_k\to0$ by \eqref{eq:grid-scales} and $|g_{j,k}-\omega_j|\leq h_k$,
$g_{j,k}\in\mathcal U_{j,k}$ for all sufficiently large $k$.

For $\omega\in\mathcal G_k\setminus
\bigcup_{j=1}^q\mathcal U_{j,k}$,
we have that $|\omega-\omega_j|\geq\rho_k$
for every $j=1,\ldots,q$. Hence \eqref{eq:kernel-bound} implies
\[
 |D_k(\omega_j-\omega)|
 \leq\frac{2}{T_k|\omega_j-\omega|}
 \leq\frac{2}{T_k\rho_k}.
\]
Applying the triangle inequality to \eqref{eq:F-decomposition},
we obtain
\[
\begin{aligned}
 \|F_k(\omega)\|
 &\leq\sum_{j=1}^q
       \|c_{j,f}\|\,|D_k(\omega_j-\omega)|
       +\|\mathcal E_k(\omega)\|\\
 &\leq\frac{2\sum_{j=1}^q\|c_{j,f}\|}{T_k\rho_k}
       +\max_{\nu\in\mathcal G_k}\|\mathcal E_k(\nu)\|.
\end{aligned}
\]
Taking the maximum over these grid points and using
\eqref{eq:uniform-remainder} therefore gives
\begin{equation}\label{eq:outside-peaks}
\begin{aligned}
 \max_{\omega\in\mathcal G_k\setminus
                    \bigcup_{j=1}^q\mathcal U_{j,k}}
       \|F_k(\omega)\|&\leq
 \frac{2\sum_{j=1}^q\|c_{j,f}\|}{T_k\rho_k}
 +o(\gamma_k)\\
 &=o(\gamma_k).
\end{aligned}
\end{equation}
where the last equality is due to \eqref{eq:scale-properties}.
The maximum in \eqref{eq:outside-peaks} is considered only
when its index set is nonempty.
Meanwhile, \eqref{eq:nearest-peak}, the nonzero coefficients
in \eqref{eq:identification-model}, and $\gamma_k\to0$ imply
\[
 \|F_k(g_{j,k})\|
 \geq\frac12\min_{1\leq\ell\leq q}\|c_{\ell,f}\|
 >\gamma_k
\]
for every $j$ and all sufficiently large $k$.

By \eqref{eq:outside-peaks} and \eqref{eq:peak-search}, all grid points outside
$\bigcup_{j=1}^q\mathcal U_{j,k}$ satisfy
$\|F_k(\omega)\|\leq\gamma_k$ for all sufficiently large $k$.
Thus every accepted grid point $\widehat\omega\in\mathcal G_k$
belongs to some $\mathcal U_{j,k}$.
Since 
$\rho_k\to0$, for large enough $k$, we have
$|\omega_\ell-\omega_j|>4\rho_k$ whenever $\ell\neq j$.

Suppose $\widehat\omega\in\mathcal U_{j,k}$ is accepted.
The deletion rule \eqref{eq:peak-deletion} removes all
remaining grid points satisfying
$|\omega-\widehat\omega|\leq2\rho_k$.
For every $\omega\in\mathcal U_{j,k}$,
\[
 |\omega-\widehat\omega|
 \leq|\omega-\omega_j|+|\omega_j-\widehat\omega|
 <2\rho_k.
\]
Hence the entire neighborhood $\mathcal U_{j,k}$ is removed,
so no further point from this neighborhood can be selected.

In contrast, for $\omega\in\mathcal U_{\ell,k}$ with $\ell\neq j$,
\[
\begin{aligned}
 |\omega-\widehat\omega|
 &\geq|\omega_\ell-\omega_j|
      -|\omega-\omega_\ell|
      -|\widehat\omega-\omega_j|\\
 &>4\rho_k-\rho_k-\rho_k=2\rho_k.
\end{aligned}
\]
Therefore, every point in the other neighborhoods remains
available for subsequent selections.

For all sufficiently large $k$, \eqref{eq:nearest-peak} and
$\gamma_k\to0$ give
$\|F_k(g_{j,k})\|>\gamma_k$ for every $j$.
If no point has yet been selected from $\mathcal U_{j,k}$,
the preceding deletion argument shows that $g_{j,k}$
is still available.
Thus the algorithm cannot stop while any neighborhood
remains unselected.

Each selection removes exactly one neighborhood and preserves
all others.
After $q$ selections, all $q$ neighborhoods have been removed.
Every remaining candidate then satisfies
$\|F_k(\omega)\|\leq\gamma_k$ by \eqref{eq:outside-peaks},
so the algorithm stops.
Therefore, $\widehat q_k=q$ for all sufficiently large $k$,
which proves \eqref{eq:order-consistency}.

Each neighborhood contributes one selected point; sorting
therefore gives $\widehat\omega_{j,k}\in\mathcal U_{j,k}$.
By the definition \eqref{eq:peak-neighborhoods},
\begin{equation}\label{eq:coarse-localization}
 |\widehat\omega_{j,k}-\omega_j|<\rho_k,
 \qquad j=1,\ldots,q,
\end{equation}
for all sufficiently large $k$.

When $q=0$, \eqref{eq:F-decomposition} gives
$F_k(\omega)=\mathcal E_k(\omega)$.
By \eqref{eq:uniform-remainder},
$\max_{\omega\in\mathcal G_k}\|F_k(\omega)\|=o(\gamma_k)$,
so no candidate exceeds the threshold for all sufficiently
large $k$.
The algorithm therefore returns $\widehat q_k=0$.

\emph{Frequency accuracy:}
Suppose $q>0$ and fix $j\in\{1,\ldots,q\}$ and
$\varepsilon>0$.
Define
\begin{equation}\label{eq:sinc-gap}
 \mu_\varepsilon
 :=\sup_{|x|\geq\varepsilon}
       \left|\frac{\sin(x/2)}{x/2}\right|<1.
\end{equation}

For $\omega\in\mathcal U_{j,k}$ and $\ell\neq j$,
\eqref{eq:peak-neighborhoods} and $\rho_k\to0$ give
\[
 |\omega_\ell-\omega|
 \geq|\omega_\ell-\omega_j|-\rho_k
 \geq\frac12|\omega_\ell-\omega_j|>0
\]
for all sufficiently large $k$.
Thus \eqref{eq:kernel-bound} gives
$D_k(\omega_\ell-\omega)=O(T_k^{-1})$ uniformly
on $\mathcal U_{j,k}$.
Together with \eqref{eq:F-decomposition} and
\eqref{eq:uniform-remainder}, this yields
\begin{equation}\label{eq:local-peak-form}
 F_k(\omega)
 =c_{j,f}D_k(\omega_j-\omega)+o(1),
 \qquad \omega\in\mathcal U_{j,k}.
\end{equation}

For $\omega\in\mathcal U_{j,k}$ satisfying
$T_k|\omega-\omega_j|\geq\varepsilon$,
\eqref{eq:kernel-bound} and \eqref{eq:sinc-gap} give
$|D_k(\omega_j-\omega)|\leq\mu_\varepsilon$.
Hence \eqref{eq:local-peak-form} and
\eqref{eq:nearest-peak} imply
\[
\begin{aligned}
 \|F_k(\omega)\|
 &\leq\mu_\varepsilon\|c_{j,f}\|+o(1),\\
 \|F_k(g_{j,k})\|
 &=\|c_{j,f}\|+o(1).
\end{aligned}
\]
Since $\mu_\varepsilon<1$ and $c_{j,f}\neq0$,
\begin{align}\label{ineq_Fkg_ome}
        \|F_k(g_{j,k})\|-\|F_k(\omega)\|
 \geq(1-\mu_\varepsilon)\|c_{j,f}\|+o(1)>0
\end{align}
uniformly over these points for all sufficiently large $k$.

When $\widehat\omega_{j,k}$ is selected,
$g_{j,k}$ belongs to the current search set $\mathcal G_k^{\mathrm{rem}}$.
Hence \eqref{eq:peak-search} gives
\[
 \|F_k(\widehat\omega_{j,k})\|
 \geq\|F_k(g_{j,k})\|.
\]
By \eqref{eq:coarse-localization},
$\widehat\omega_{j,k}\in\mathcal U_{j,k}$.
If $T_k|\widehat\omega_{j,k}-\omega_j|\geq\varepsilon$,
then \eqref{ineq_Fkg_ome}, with
$\omega=\widehat\omega_{j,k}$, gives
\[
 \|F_k(\widehat\omega_{j,k})\|
 <\|F_k(g_{j,k})\|,
\]
contradicting the maximization rule \eqref{eq:peak-search}.
Therefore,
\(
 T_k|\widehat\omega_{j,k}-\omega_j|<\varepsilon
\)
for all sufficiently large $k$.
Since $\varepsilon>0$ is arbitrary, this proves
\eqref{eq:frequency-consistency}.
The argument holds on the probability-one event fixed above,
so the convergence is almost sure.
\end{IEEEproof}



Theorem~\ref{thm:identification} supplies exactly the information
needed for the next step: the filter length can be chosen from
$\widehat q_k$, and its output gain from the sorted frequency
tuple. The amplitudes used to detect peaks do not enter either
control calculation.

\begin{remark}
\label{rem:window-tradeoff}
The growing blocks in \eqref{eq:blocks} provide the noise
averaging used in Theorem~\ref{thm:identification}, but also
increase the time between successive controller updates.
For practical use, the observation length may be capped at
$T_{\max}\geq t_{\mathrm b}$ by setting
\[
\begin{aligned}
 \widetilde T_k&=\min\{t_{\mathrm b}2^k,T_{\max}\},\\
 \widetilde t_{k+1}&=\widetilde t_k+\widetilde T_k,
 \qquad \widetilde t_0=t_{\mathrm b}.
\end{aligned}
\]
The observation blocks then become
$[\widetilde t_k,\widetilde t_{k+1}]$.
Once the cap is reached, the search grid, detection threshold,
and exclusion radius should be held fixed at values selected
for $T_{\max}$.
The fixed grid keeps the number of Fourier accumulators
constant, so the storage requirement and computational cost
per integration step no longer grow with time.
A smaller $T_{\max}$ permits more frequent updates and may
improve the response to frequency jumps, at the cost of reduced
noise averaging and frequency resolution.
The cap and detection parameters should therefore be chosen
according to the required identification accuracy and
available computational resources.
\end{remark}

\section{Regulation with Known Frequencies}
\label{sec:controller}
We first suppose that the active positive frequencies
$\Omega_d=\{\omega_1,\ldots,\omega_q\}$ in
\eqref{eq:harmonic-d} are known, with
$0<\omega_1<\cdots<\omega_q$ and $q=0$ allowed.
Their amplitudes and phases remain unknown.
We construct a fixed error-feedback controller and establish its
stability and regulation properties. Section~\ref{sec:closed-loop}
then replaces these frequencies by their online estimates.

\subsection{Controller realization and closed-loop equations}
Under Assumption~\ref{ass:plant-exosystem}, choose
$K\in\mathbb R^{m\times n}$ such that
\begin{equation}\label{eq:AK}
 A_K:=A+BK\ \text{is Hurwitz},\qquad C_K:=C+DK.
\end{equation}
Retain $L_o$ and $\alpha>0$ from \eqref{eq:Ao} and
\eqref{eq:lowpass-implementation}, and choose $0<\kappa<1$.
For each integer $r\geq1$, define
\begin{equation}\label{eq:filter-shift}
\begin{gathered}
 N_r=
 \begin{bmatrix}
 0_{1\times(r-1)}&0\\
 I_{r-1}&0_{(r-1)\times1}
 \end{bmatrix}\quad(r\geq2),\qquad N_1=[0],\\
 \mathbf e_1^{(r)}=\col(1,0,\ldots,0)\in\mathbb R^r.
\end{gathered}
\end{equation}
Here $r$ depends on the number of frequencies in $d$
and will be specified in the next subsection.

Set
\begin{equation}\label{eq:variable-filter-matrices}
 \mathsf A_{f,r}=\alpha(\kappa N_r-I_r)\otimes I_p,\qquad
 \mathsf B_{f,r}=\alpha\mathbf e_1^{(r)}\otimes I_p.
\end{equation}
For $\eta=\col(\eta_1,\ldots,\eta_r)\in\mathbb R^{pr}$,
the cascade driven by the residual $\mathsf R$ in
\eqref{eq:integrated-residual} is
\begin{equation}\label{eq:controller-filter}
\begin{aligned}
 \dd\eta_1&=-\alpha\eta_1\,\dd t+\alpha\,\dd\mathsf R,\\
 \dot\eta_\ell&=-\alpha\eta_\ell+\alpha\kappa\eta_{\ell-1},
                       \qquad 2\leq\ell\leq r,
\end{aligned}
\end{equation}
or, equivalently,
\begin{equation}\label{eq:filter-realization}
 \dd\eta=\mathsf A_{f,r}\eta\,\dd t
                         +\mathsf B_{f,r}\,\dd\mathsf R.
\end{equation}
Taking $\eta_1(0)=y_f(0)$ gives $\eta_1=y_f$ by
\eqref{eq:lowpass-implementation}.
Since $\|N_r\|\leq1$, we have $N_r+N_r^\top\preceq2I_r$.
The Kronecker product with $I_p$ preserves this inequality; hence
\begin{equation}\label{eq:uniform-filter-dissipation}
\begin{aligned}
 \mathsf A_{f,r}+\mathsf A_{f,r}^\top
 &=\alpha[\kappa(N_r+N_r^\top)-2I_r]\otimes I_p\\
 &\preceq-2bI_{pr},\qquad b:=\alpha(1-\kappa)>0,\\
 \|e^{\mathsf A_{f,r}t}\|&\leq e^{-bt},\qquad t\geq0.
\end{aligned}
\end{equation}
For $\dot x=\mathsf A_{f,r}x$, the first inequality gives
\[
 \frac{\dd}{\dd t}\|x(t)\|^2\leq-2b\|x(t)\|^2.
\]
Integration yields the exponential estimate.
The constant $b$ is independent of the filter order.

For a fixed $H=[H_1\ \cdots\ H_r]\in\mathbb R^{m\times pr}$,
combine \eqref{eq:implementable-observer},
\eqref{eq:integrated-residual}, and \eqref{eq:filter-realization}
to obtain the controller
\begin{equation}\label{eq:complete-controller}
\begin{aligned}
 \dd\xi&=[A_o\xi+(B-L_oD)u]\,\dd t+L_o\,\dd E,\\
 \dd\eta&=[\mathsf A_{f,r}\eta
            -\mathsf B_{f,r}(C\xi+Du)]\,\dd t
            +\mathsf B_{f,r}\,\dd E,\\
 u&=K\xi+H\eta,
\end{aligned}
\end{equation}
where the gain $H$ will be determined from \eqref{eq:gain-interpolation} in the next subsection.
Its input is the measured process $E$ and its output is $u$.
Substituting the last equation into the first two yields
\begin{equation}\label{eq:controller-state-realization}
\begin{aligned}
 \dd x_{\mathrm{ctr}}&=A_c(r,H)x_{\mathrm{ctr}}\,\dd t+B_c(r)\,\dd E,\\
 u&=K_c(H)x_{\mathrm{ctr}},\qquad x_{\mathrm{ctr}}=\col(\xi,\eta),
\end{aligned}
\end{equation}
where
\begin{equation}\label{eq:controller-state-matrices}
\begin{aligned}
 A_c(r,H)&=
 \begin{bmatrix}
 A_K-L_oC_K&(B-L_oD)H\\
 -\mathsf B_{f,r}C_K&\mathsf A_{f,r}-\mathsf B_{f,r}DH
 \end{bmatrix},\\
 B_c(r)&=\begin{bmatrix}L_o\\\mathsf B_{f,r}\end{bmatrix},
 \qquad K_c(H)=\begin{bmatrix}K&H\end{bmatrix}.
\end{aligned}
\end{equation}
Thus the controller has no direct feedthrough from $\dd E$ to $u$.

We next connect \eqref{eq:complete-controller} to the equivalent
plant \eqref{eq:equivalent-plant}. Recall
$\bar z=z-\Pi w$ and $\chi=\bar z-\xi$ from
\eqref{eq:equivalent-definitions} and \eqref{eq:chi}.
Using \eqref{eq:residual-system}, we obtain the closed-loop system:
\begin{equation}\label{eq:fixed-closed-loop}
\begin{aligned}
 \dd\chi&=A_o\chi\,\dd t-L_oV_m\,\dd W_m+V_p\,\dd W_p,\\
 \dd\eta&=(\mathsf A_{f,r}\eta+\mathsf B_{f,r}C\chi
                   +\mathsf B_{f,r}d)\,\dd t
                   +\mathsf B_{f,r}V_m\,\dd W_m,\\
 \dd\bar z&=(A_K\bar z-BK\chi+BH\eta+L_od)\,\dd t
                   +V_p\,\dd W_p,\\
 e_0&=C_K\bar z-DK\chi+DH\eta+d.
\end{aligned}
\end{equation}
Define $x_{\mathrm{cl}}=\col(\chi,\eta,\bar z)$, 
\begin{equation}\label{eq:closed-loop-matrices}
\begin{gathered}
 \mathsf A_{r,H}=
 \begin{bmatrix}
 A_o&0&0\\
 \mathsf B_{f,r}C&\mathsf A_{f,r}&0\\
 -BK&BH&A_K
 \end{bmatrix},\qquad
 \mathsf D_r=\begin{bmatrix}0\\\mathsf B_{f,r}\\L_o\end{bmatrix},\\
 \mathsf B_{m,r}=\begin{bmatrix}-L_o\\\mathsf B_{f,r}\\0\end{bmatrix},\quad
 \mathsf B_{p,r}=\begin{bmatrix}I_n\\0_{pr\times n}\\I_n\end{bmatrix},\\
 \mathsf C_{r,H}=\begin{bmatrix}-DK&DH&C_K\end{bmatrix}.
\end{gathered}
\end{equation}
Then \eqref{eq:fixed-closed-loop} becomes
\begin{equation}\label{eq:closed-loop-state-form}
\begin{aligned}
 \dd x_{\mathrm{cl}}&=(\mathsf A_{r,H}x_{\mathrm{cl}}+\mathsf D_rd)\,\dd t
                             +\mathsf B_{m,r}V_m\,\dd W_m
                             +\mathsf B_{p,r}V_p\,\dd W_p,\\
 e_0&=\mathsf C_{r,H}x_{\mathrm{cl}}+d.
\end{aligned}
\end{equation}
Clearly,
\begin{equation}\label{eq:fixed-closed-loop-spectrum}
 \sigma(\mathsf A_{r,H})
 =\sigma(A_o)\cup\sigma(\mathsf A_{f,r})\cup\sigma(A_K)
 \subset\{s:\Rea s<0\}.
\end{equation}
Thus every fixed finite pair $(r,H)$ gives an internally stable
closed loop. We next choose $(r,H)$ to eliminate the steady-state
tracking error caused by $d$.

\subsection{Frequency cancellation and gain construction}
To calculate the transfer from $d$ to $e_0$, set $V_m=V_p=0$ and
take zero initial states in \eqref{eq:fixed-closed-loop}.
Equation~\eqref{eq:residual-system} gives $\chi=0$ and
$\dd\mathsf R=d\,\dd t$. The resulting closed-loop equations are
\begin{equation}\label{eq:zero-state-closed-loop}
\begin{aligned}
 \dot{\bar z}&=A_K\bar z+BH\eta+L_od,\\
 \dot\eta&=\mathsf A_{f,r}\eta+\mathsf B_{f,r}d,\\
 e_0&=C_K\bar z+DH\eta+d.
\end{aligned}
\end{equation}
Here $d$ is treated as the input and $e_0$ as the output.
Let $v_c:=H\eta$. The first and third
equations give
\begin{equation}\label{eq:stabilized-transfers}
\begin{aligned}
 e_0(s)&=\mathscr P_K(s)v_c(s)+\mathscr D_K(s)d(s),\\
 \mathscr P_K(s)&=C_K(sI-A_K)^{-1}B+D,\\
 \mathscr D_K(s)&=C_K(sI-A_K)^{-1}L_o+I_p.
\end{aligned}
\end{equation}
To calculate the remaining transfer, use
$h_f(s)=\alpha/(s+\alpha)$ from \eqref{eq:lowpass-transfer}.
The cascade \eqref{eq:controller-filter} with
$\dd\mathsf R=d\,\dd t$ and zero initial states gives
\[
 \eta_1(s)=h_f(s)d(s),\qquad
 \eta_\ell(s)=\kappa h_f(s)\eta_{\ell-1}(s).
\]
Consequently, $\eta_\ell(s)=\kappa^{\ell-1}h_f(s)^\ell d(s)$.
The transfer matrix from $d$ to $v_c$ in
\eqref{eq:zero-state-closed-loop} is therefore
\begin{equation}\label{eq:stable-interpolant}
\begin{aligned}
 v_c(s)&=\mathscr F_H(s)d(s),\\
 \mathscr F_H(s)
 &=H(sI-\mathsf A_{f,r})^{-1}\mathsf B_{f,r}=\sum_{\ell=1}^{r}H_\ell\kappa^{\ell-1}h_f(s)^\ell.
\end{aligned}
\end{equation}
Substitution into \eqref{eq:stabilized-transfers} gives the
closed-loop transfer matrix of \eqref{eq:closed-loop-state-form} from $d$ to $e_0$:
\begin{equation}\label{eq:regulation-transfer}
\begin{aligned}
 e_0(s)&=\mathscr T_H(s)d(s),\\
 \mathscr T_H(s)&=\mathscr D_K(s)+\mathscr P_K(s)\mathscr F_H(s).
\end{aligned}
\end{equation}
Since the closed loop is stable by
\eqref{eq:fixed-closed-loop-spectrum}, to have a zero steady-state tracking error, it suffices to let $\mathscr T_H(\ii\omega)=0$ for $\omega$ being the frequencies of $d(t)$,  see for instance, \cite{Bymes2000},  \cite{Francis1977}.

To account for the possible constant component in the external signal, define
\begin{equation}\label{eq:zero-node-selector}
 \iota_0=
 \begin{cases}
 1,&\rank\mathcal R_u(0)=n+p,\\
 0,&\rank\mathcal R_u(0)<n+p.
 \end{cases}
\end{equation}
The rank condition in Assumption~\ref{ass:regulation} is imposed
only at $\lambda\in\sigma(S)$. It implies
\begin{equation}\label{eq:zero-node-coverage}
 \iota_0=0\ \Longrightarrow\ 0\notin\sigma(S)
                         \ \Longrightarrow\ d_0=0,
\end{equation}
where $d_0$ is the constant component in \eqref{eq:harmonic-d}.

Let $0<\omega_1<\cdots<\omega_q$ be the known active positive
frequencies in \eqref{eq:harmonic-d}, with $q=0$ allowed. Set
\begin{equation}\label{eq:interpolation-nodes}
\begin{aligned}
 \Lambda&=
 \begin{cases}
 \{0,\pm\ii\omega_1,\ldots,\pm\ii\omega_q\},
                                           &\iota_0=1,\\
 \{\pm\ii\omega_1,\ldots,\pm\ii\omega_q\},
                                           &\iota_0=0,
 \end{cases}\\
 r&=\max\{1,|\Lambda|\}
              =\max\{1,2q+\iota_0\}.
\end{aligned}
\end{equation}
If $\iota_0=q=0$, the node set is empty and the retained
stage $r=1$ supplies only $\eta_1=y_f$.

The state-feedback identity
\[
 \mathcal R_u(\lambda)
 \begin{bmatrix}I_n&0\\K&I_m\end{bmatrix}
 =\begin{bmatrix}\lambda I-A_K&-B\\C_K&D\end{bmatrix}
\]
and elimination of the invertible block $\lambda I-A_K$ give
the standard rank relation \cite{Trentelman2001}
\begin{equation}\label{eq:rank-transfer-equivalence}
 \rank\mathcal R_u(\lambda)=n+\rank\mathscr P_K(\lambda),
 \qquad\lambda\in\ii\mathbb R.
\end{equation}
Since $\pm\ii\omega_j\in\sigma(S)$,
Assumption~\ref{ass:regulation}, \eqref{eq:zero-node-selector},
and \eqref{eq:rank-transfer-equivalence} imply
\begin{equation}\label{eq:known-node-rank}
 \rank\mathscr P_K(\lambda)=p,\qquad\lambda\in\Lambda.
\end{equation}
For $\Lambda\neq\varnothing$, choose
\begin{equation}\label{eq:interpolation-target}
 \mathscr M(\lambda)
 =-\mathscr P_K(\lambda)^*
 [\mathscr P_K(\lambda)\mathscr P_K(\lambda)^*]^{-1}
 \mathscr D_K(\lambda).
\end{equation} 
It satisfies $\mathscr P_K(\lambda)\mathscr M(\lambda)
=-\mathscr D_K(\lambda)$. The gain
$H=[H_1\ \cdots\ H_{r}]$ is determined by
\begin{equation}\label{eq:gain-interpolation}
 \sum_{\ell=1}^{r}
 H_\ell\kappa^{\ell-1}h_f(\lambda)^\ell
 =\mathscr M(\lambda),\qquad\lambda\in\Lambda.
\end{equation}
That is, writing
$\Lambda=\{\lambda_1,\ldots,\lambda_r\}\subset\ii\mathbb R$,
\begingroup
\footnotesize
\setlength{\arraycolsep}{3pt}
\renewcommand{\arraystretch}{0.9}

\[
 (\mathcal V_r\otimes I_m)
 \begin{bmatrix}
 H_1\\ \vdots\\ H_{r}
 \end{bmatrix}
 =
 \begin{bmatrix}
 \mathscr M(\lambda_1)\\
 \vdots\\
 \mathscr M(\lambda_{r})
 \end{bmatrix},
\]\endgroup
where
\[
 \mathcal V_r=
 \begin{bmatrix}
 h_f(\lambda_1)
 &\kappa h_f(\lambda_1)^2
 &\cdots
 &\kappa^{r-1}h_f(\lambda_1)^{r}\\
 \vdots&\vdots&\ddots&\vdots\\
 h_f(\lambda_{r})
 &\kappa h_f(\lambda_{r})^2
 &\cdots
 &\kappa^{r-1}
 h_f(\lambda_{r})^{r}
 \end{bmatrix}.
\]

Noting $h_f(s)=\alpha/(s+\alpha)$ and $\lambda_i\neq \lambda_j $, for $i\neq j$, we have
 $h_f(\lambda_i)\neq h_f(\lambda_j)$. The Vandermonde formula gives
\begin{equation}\label{eq:interpolation-Vandermonde}
\resizebox{0.9\columnwidth}{!}{$\displaystyle
 \det\mathcal V_r
 =
 \kappa^{r(r-1)/2}
 \prod_{j=1}^{r}h_f(\lambda_j)
 \prod_{1\leq j<i\leq r}
 \bigl[h_f(\lambda_i)-h_f(\lambda_j)\bigr]
 \ne0.
$}
\end{equation}
Thus $\mathcal V_r$ is invertible, and
\eqref{eq:gain-interpolation} uniquely determines
$H_1,\ldots,H_{r}$. 
By conjugate symmetry, $\overline{H}$ also satisfies
\eqref{eq:gain-interpolation}; uniqueness therefore implies
that $H$ is real.
For fixed $q$, $H$ depends continuously on
$(\omega_1,\ldots,\omega_q)$ while the nodes remain distinct
and \eqref{eq:known-node-rank} holds: both the targets in
\eqref{eq:interpolation-target} and $\mathcal V_r^{-1}$
are continuous.

If $\Lambda=\varnothing$, then $r=1$,
and we set
\begin{equation}\label{eq:empty-interpolation-gain}
 H=0_{m\times p}.
\end{equation}

Using \eqref{eq:stable-interpolant},
\eqref{eq:interpolation-target}, and \eqref{eq:gain-interpolation}
now gives
\begin{equation}\label{eq:interpolation-cancellation}
\begin{aligned}
 \mathscr T_H(\lambda)
 &=\mathscr D_K(\lambda)
       +\mathscr P_K(\lambda)\mathscr M(\lambda)=0_{p\times p},\ \  \lambda\in\Lambda.
\end{aligned}
\end{equation}

The \textit{known-frequency controller} is \eqref{eq:complete-controller}
with $r$ from \eqref{eq:interpolation-nodes} and $H$ from
\eqref{eq:gain-interpolation} or \eqref{eq:empty-interpolation-gain}.
If $q=\iota_0=0$, then \eqref{eq:zero-node-coverage} gives
$d=0$, and the construction uses $r=1$, $H=0$.

\begin{remark}[Internal model in the controller]
\label{rem:controller-internal-model}
Suppose $\Lambda\neq\varnothing$, and let $r,H$ be given by
\eqref{eq:interpolation-nodes} and \eqref{eq:gain-interpolation}.
We show that the controller matrix in
\eqref{eq:controller-state-matrices} satisfies
\begin{equation}\label{eq:controller-internal-model}
 \dim\ker(\lambda I-A_c(r,H))\geq p,
 \qquad \lambda\in\Lambda.
\end{equation}
Abbreviate $A_f=\mathsf A_{f,r}$ and $B_f=\mathsf B_{f,r}$,
and define
\[
 A_b=\begin{bmatrix}A_K&BH\\0&A_f\end{bmatrix},
 \qquad C_b=\begin{bmatrix}C_K&DH\end{bmatrix}.
\]
Using $B_c(r)=\col(L_o,B_f)$, block inversion gives
\begin{equation} \label{sssp}
	\begin{aligned}
 &I_p+C_b(sI-A_b)^{-1}B_c(r)=I_p+C_K(sI-A_K)^{-1}L_o\\
 &\qquad+
 [C_K(sI-A_K)^{-1}B+D]H(sI-A_f)^{-1}B_f\\
 &=\mathscr D_K(s)+\mathscr P_K(s)\mathscr F_H(s)=\mathscr T_H(s),
\end{aligned}
\end{equation}
where the last two equalities follow from
\eqref{eq:stabilized-transfers},
\eqref{eq:stable-interpolant}, and
\eqref{eq:regulation-transfer}.
From the expression for $A_c(r,H)$ in
\eqref{eq:controller-state-matrices}, we have
\begin{equation}\label{eq:internal-model-factorization}
\begin{aligned}
 A_c(r,H)&=A_b-B_c(r)C_b.
\end{aligned}
\end{equation}
For $\lambda\in\Lambda\subset\ii\mathbb R$,
 the matrix $\lambda I-A_b$ is invertible, since $A_b$ is Hurwitz.
Set
 $Z_\lambda=(\lambda I-A_b)^{-1}B_c(r).$ Using
\eqref{sssp}--\eqref{eq:internal-model-factorization},
together with \eqref{eq:interpolation-cancellation}, we have
\begin{equation}
	\begin{aligned}
 A_bZ_\lambda&=\lambda Z_\lambda-B_c(r),\ 
 C_bZ_\lambda=\mathscr T_H(\lambda)-I_p=-I_p.
\end{aligned}
\end{equation}
Then
\[
\begin{aligned}
 A_c(r,H)Z_\lambda
 &=A_bZ_\lambda-B_c(r)C_bZ_\lambda=\lambda Z_\lambda-B_c(r)+B_c(r)\\
 &=\lambda Z_\lambda, \quad \lambda\in \Lambda.
\end{aligned}
\]
Thus every column of $Z_\lambda$ belongs to
$\ker(\lambda I-A_c(r,H))$.
Moreover, $C_bZ_\lambda=-I_p$ implies
\[
 p=\rank(C_bZ_\lambda)\leq\rank Z_\lambda\leq p,
\]
which proves \eqref{eq:controller-internal-model}.
This is the classical
$p$-copy internal-model condition \cite{FrancisWonham1976} for robust output regulation.

\end{remark}

\subsection{Stability and noise-dependent regulation performance}
Use the subscript $\star$ for this fixed known-frequency
closed loop, which will serve as the comparison system in
Section~\ref{sec:closed-loop}. With $r,H$ defined above, set
\begin{equation}\label{eq:limiting-matrices}
\begin{gathered}
 \mathsf A_\star=\mathsf A_{r,H},\quad
 \mathsf B_{m,\star}=\mathsf B_{m,r},\quad
 \mathsf B_{p,\star}=\mathsf B_{p,r},\\
 \mathsf C_\star=\mathsf C_{r,H},\qquad
 \mathsf D_\star=\mathsf D_r.
\end{gathered}
\end{equation}
Its state and true error are
$x_{\mathrm{cl},\star}=\col(\chi_\star,\eta_\star,\bar z_\star)$ and
$e_{0,\star}=\mathsf C_\star x_{\mathrm{cl},\star}+d$.
For a fixed active frequency set, the matrices in
\eqref{eq:limiting-matrices} are independent of the harmonic
amplitudes and both noise covariances. Define
\begin{equation}\label{eq:closed-loop-noise-intensity}
 \mathsf N_\star=
 \mathsf B_{m,\star}\Sigma_m\mathsf B_{m,\star}^\top
 +\mathsf B_{p,\star}\Sigma_p\mathsf B_{p,\star}^\top.
\end{equation}
The two nonzero blocks of $\mathsf B_{p,\star}$ are driven by
the same $W_p$, so this expression includes their cross covariance.

Since $\mathsf A_\star$ is Hurwitz, the standard linear
stochastic-system theory \cite[Sec.~6.5]{SarkkaSolin2019}
gives the stationary state-noise covariance
\begin{equation}\label{eq:stationary-covariance-integral}
 \Sigma_{\mathrm{cl},\star}=
 \int_0^\infty e^{\mathsf A_\star \tau}\mathsf N_\star
 e^{\mathsf A_\star^\top \tau}\,\dd \tau.
\end{equation}
It is the unique positive semidefinite solution of
\begin{equation}\label{eq:closed-loop-Lyapunov}
 \mathsf A_\star\Sigma_{\mathrm{cl},\star}+
 \Sigma_{\mathrm{cl},\star}\mathsf A_\star^\top+\mathsf N_\star=0.
\end{equation}
Define
\begin{equation}\label{eq:limiting-error-covariance}
 \Sigma_e=\mathsf C_\star\Sigma_{\mathrm{cl},\star}\mathsf C_\star^\top,
 \qquad J=\operatorname{tr}(\Sigma_e),
\end{equation}
and the deterministic matrices
\begin{equation}\label{eq:error-noise-intensity}
\begin{aligned}
 \mathcal K_{e,m}(\tau)&=\mathsf C_\star e^{\mathsf A_\star \tau}
                                      \mathsf B_{m,\star},\quad
 \mathcal K_{e,p}(\tau)=\mathsf C_\star e^{\mathsf A_\star \tau}
                                      \mathsf B_{p,\star},\\
 \mathcal M_{e,m}&=\int_0^\infty
                  \mathcal K_{e,m}(\tau)^\top\mathcal K_{e,m}(\tau)\,\dd \tau
                  \succeq0,\\
 \mathcal M_{e,p}&=\int_0^\infty
                  \mathcal K_{e,p}(\tau)^\top\mathcal K_{e,p}(\tau)\,\dd \tau
                  \succeq0.
\end{aligned}
\end{equation}
Equations \eqref{eq:stationary-covariance-integral} and
\eqref{eq:limiting-error-covariance}, together with independence
of $W_m,W_p$ and the cyclic property of the trace, yield
\begin{equation}\label{eq:noise-performance-functional}
\begin{aligned}
 \Sigma_e&=\int_0^\infty\bigl[
       \mathcal K_{e,m}(\tau)\Sigma_m\mathcal K_{e,m}(\tau)^\top\\
 &\qquad{}+\mathcal K_{e,p}(\tau)\Sigma_p
                   \mathcal K_{e,p}(\tau)^\top\bigr]\,\dd \tau,\\
 J_m&=\operatorname{tr}(\mathcal M_{e,m}\Sigma_m),\qquad
 J_p=\operatorname{tr}(\mathcal M_{e,p}\Sigma_p),\\
 J&=J_m+J_p\leq\|\mathcal M_{e,m}\|\operatorname{tr}(\Sigma_m)
       +\|\mathcal M_{e,p}\|\operatorname{tr}(\Sigma_p).
\end{aligned}
\end{equation}
Thus $J$ is a nonnegative linear function of $(\Sigma_m,\Sigma_p)$.
Both covariances are used only to evaluate performance; neither
is required by the controller or frequency estimator.

\begin{theorem}[Regulation with known frequencies]
\label{thm:known-frequency-regulation}
Suppose Assumptions~\ref{ass:plant-exosystem} and
\ref{ass:regulation} hold.
Let $K$ and $L_o$ be such that $A_K$ and $A_o$ are Hurwitz,
and fix $\alpha>0$ and $0<\kappa<1$.
Consider the controller \eqref{eq:complete-controller},
with the order specified by \eqref{eq:interpolation-nodes}
and the gain determined by \eqref{eq:gain-interpolation}
when $\Lambda\neq\varnothing$, or by
\eqref{eq:empty-interpolation-gain} otherwise.

For any deterministic filter initial state $\eta(0)$,
the closed loop \eqref{eq:closed-loop-state-form} has a unique global strong solution and, for $a=2,4$,
\begin{equation}\label{eq:known-frequency-moment-bound}
\begin{split}
 \sup_{t\geq0}\E\bigl[
 &\|x_{\mathrm{cl},\star}(t)\|^a+\|z_\star(t)\|^a
   +\|\xi_\star(t)\|^a\\
 &+\|u_\star(t)\|^a+\|e_{0,\star}(t)\|^a\bigr]<\infty.
\end{split}
\end{equation}
Moreover, as $t\to\infty$,
\begin{equation}\label{eq:known-frequency-regulation}
\begin{aligned}
 \E[e_{0,\star}(t)]&\longrightarrow0,\quad
 \E[e_{0,\star}(t)e_{0,\star}(t)^\top]&\longrightarrow \Sigma_e,
\end{aligned}
\end{equation}
and
\begin{equation}\label{eq:known-frequency-mse}
\begin{aligned}
 \lim_{t\to\infty}\E\|e_{0,\star}(t)\|^2
 &=J_m+J_p=J\\
 \leq &\|\mathcal M_{e,m}\|\operatorname{tr}(\Sigma_m)
       +\|\mathcal M_{e,p}\|\operatorname{tr}(\Sigma_p),
\end{aligned}
\end{equation}
where $\Sigma_e,J$ and $J_m,J_p$ are defined by
\eqref{eq:limiting-error-covariance} and \eqref{eq:noise-performance-functional}.
If $V_m=V_p=0$, then $e_{0,\star}(t)$ decays exponentially.
\end{theorem}
\begin{IEEEproof}
Consider  the closed-loop system
\eqref{eq:closed-loop-state-form} with the matrices
defined in \eqref{eq:limiting-matrices}.
By \eqref{eq:fixed-closed-loop-spectrum}, there are
$M_\star,\beta_\star>0$ such that
$\|e^{\mathsf A_\star t}\|\leq M_\star e^{-\beta_\star t}$.
Since $d$ in \eqref{eq:harmonic-d} is bounded, the drift and
diffusion coefficients of \eqref{eq:closed-loop-state-form}
are globally Lipschitz in the state. Hence the SDE has a unique global strong
solution; see \cite[Sec.~4.5]{SarkkaSolin2019}.

\emph{Cancellation of the deterministic harmonics:}
By \eqref{eq:harmonic-d}, \eqref{eq:zero-node-coverage}, and
\eqref{eq:interpolation-nodes}, write
\(
 d(t)=\sum_{\lambda\in\Lambda}d_\lambda e^{\lambda t}.
\)
Define the particular solution of \eqref{eq:closed-loop-state-form}
\begin{equation}\label{eq:ideal-harmonic-state}
 x_d(t)=\sum_{\lambda\in\Lambda}
   (\lambda I-\mathsf A_\star)^{-1}\mathsf D_\star
                       d_\lambda e^{\lambda t}.
\end{equation}
Every resolvent exists because $\mathsf A_\star$ is Hurwitz.
Differentiating term by term gives
$\dot x_d=\mathsf A_\star x_d+\mathsf D_\star d$.
From \eqref{eq:closed-loop-state-form} and \eqref{eq:regulation-transfer}, we know
\begin{equation}\label{eq:ideal-regulation-transfer}
 I_p+\mathsf C_\star(sI-\mathsf A_\star)^{-1}\mathsf D_\star
 =\mathscr T_H(s).
\end{equation}
From \eqref{eq:interpolation-cancellation},  we obtain
\[
 \mathsf C_\star x_d(t)+d(t)
 =\sum_{\lambda\in\Lambda}
        \mathscr T_H(\lambda)d_\lambda e^{\lambda t}
 =0.
\]

\emph{Noise response and moment limits:}
Applying the variation-of-constants formula to  \eqref{eq:closed-loop-state-form} gives
\begin{equation}\label{eq:known-frequency-solution}
\begin{aligned}
 x_{\mathrm{cl},\star}(t)
 &=x_d(t)+e^{\mathsf A_\star t}[x_{\mathrm{cl},\star}(0)-x_d(0)]
                    +\zeta_\star(t),\\
 \zeta_\star(t)
 &=\int_0^t e^{\mathsf A_\star(t-s)}
                   \mathsf B_{m,\star}V_m\,\dd W_m(s)\\
 &\quad+\int_0^t e^{\mathsf A_\star(t-s)}
                   \mathsf B_{p,\star}V_p\,\dd W_p(s).
\end{aligned}
\end{equation}
Here $\zeta_\star(t)$ is a centered Gaussian vector. Independence of $W_m,W_p$, the It\^o isometry, and a change of
variables give its covariance
\begin{equation}\label{eq:known-frequency-noise-covariance}
\begin{aligned}
 \Sigma_{\mathrm{cl},\star}(t)
 &:=\E[\zeta_\star(t)\zeta_\star(t)^\top]\\
 &=\int_0^t e^{\mathsf A_\star \tau}\mathsf N_\star
       e^{\mathsf A_\star^\top \tau}\,\dd \tau
 \longrightarrow\Sigma_{\mathrm{cl},\star}.
\end{aligned}
\end{equation}
This integral is bounded uniformly in $t$.
Since $\zeta_\star(t)$ is centered Gaussian with covariance
$\Sigma_{\mathrm{cl},\star}(t)$, the fourth-moment identity (see the third equation in
\cite{Isserlis1918}, \cite[Sec.~8.2.4]{PetersenPedersen2012}) implies 
\[
 \E\|\zeta_\star(t)\|^4
 =[\operatorname{tr}\Sigma_{\mathrm{cl},\star}(t)]^2
                       +2\operatorname{tr}(\Sigma_{\mathrm{cl},\star}(t)^2)
\]
This expression is uniformly bounded. The deterministic terms in
\eqref{eq:known-frequency-solution} are bounded, so the second
and fourth moments of $x_{\mathrm{cl},\star}$ are bounded.
Noting $\xi_\star=\bar z_\star-\chi_\star$,
$u_\star=K\xi_\star+H\eta_\star$, and
$z_\star=\bar z_\star+\Pi w$, we have
\eqref{eq:known-frequency-moment-bound} holds.

Substituting \eqref{eq:known-frequency-solution} into the
output equation in \eqref{eq:closed-loop-state-form} and
using $\mathsf C_\star x_d(t)+d(t)=0$, we obtain
\begin{equation}\label{eq:ideal-error-decomposition}
\begin{aligned}
 e_{0,\star}(t)&=m_e(t)+\mathsf C_\star\zeta_\star(t),\\
 m_e(t)&:=\mathsf C_\star e^{\mathsf A_\star t}
                                  [x_{\mathrm{cl},\star}(0)-x_d(0)].
\end{aligned}
\end{equation}
Since $m_e(t)\to0$ exponentially as $t\to \infty$, we have $\E[e_{0,\star}(t)]=m_e(t)\to 0,$ and 
\[
\begin{aligned}
 \E[e_{0,\star}(t)e_{0,\star}(t)^\top]
 &=m_e(t)m_e(t)^\top+\mathsf C_\star \Sigma_{\mathrm{cl},\star}(t)\mathsf C_\star^\top
 \to  \Sigma_e.
\end{aligned}
\]
This proves \eqref{eq:known-frequency-regulation}.
Taking the trace and using \eqref{eq:noise-performance-functional}
gives \eqref{eq:known-frequency-mse}.
\end{IEEEproof}

\section{Adaptive Regulation with Unknown Frequencies}
\label{sec:closed-loop}
We now use the estimates from Theorem~\ref{thm:identification}
to update the order and gain of \eqref{eq:complete-controller}.
The fixed known-frequency system in
Theorem~\ref{thm:known-frequency-regulation} will serve as a
comparison system. Denote its order and gain by
\begin{equation}\label{eq:ideal-gain}
 r_\star=\max\{1,2q+\iota_0\},\qquad H_\star=H,
\end{equation}
where $H$ is constructed at the true frequencies
$\omega_1,\ldots,\omega_q$ by \eqref{eq:gain-interpolation}
or \eqref{eq:empty-interpolation-gain}.
These quantities are used only for comparison.
We first establish bounds during all updates, then prove
convergence to this system in mean square.

\subsection{Online order and gain updates}\label{subsec:adaptive-updates}
At $t_{k+1}$, the observation block $I_k$ in \eqref{eq:blocks}
supplies $\widehat q_k$ and
$\widehat\omega_{1,k}<\cdots<\widehat\omega_{\widehat q_k,k}$.
Using $\iota_0$ from \eqref{eq:zero-node-selector}, define
\begin{equation}\label{eq:estimated-interpolation-nodes}
 \widehat\Lambda_k=
 \begin{cases}
 \{0,\pm\ii\widehat\omega_{1,k},\ldots,
       \pm\ii\widehat\omega_{\widehat q_k,k}\},&\iota_0=1,\\
 \{\pm\ii\widehat\omega_{1,k},\ldots,
       \pm\ii\widehat\omega_{\widehat q_k,k}\},&\iota_0=0.
 \end{cases}
\end{equation}
The filter order $r$ in \eqref{eq:controller-filter} is updated by
\begin{equation}\label{eq:online-filter-order}
\begin{aligned}
 r_k&=\max\{1,|\widehat\Lambda_k|\}
      =\max\{1,2\widehat q_k+\iota_0\}, \quad
 r_{-1}=1,\\
 r(t)&=
 \begin{cases}
 1,&0\leq t<t_1,\\
 r_k,&t_{k+1}\leq t<t_{k+2},\quad k\geq0.
 \end{cases}
\end{aligned}
\end{equation}
Initialize $\eta_1(0)=y_f(0)$. At each update, retain the common
filter states, discard excess states, and set newly added states
to zero:
\begin{equation}\label{eq:filter-resizing}
 \eta_\ell(t_{k+1}^+)=
 \begin{cases}
 \eta_\ell(t_{k+1}^-),&1\leq \ell\leq\min(r_{k-1},r_k),\\
 0,&r_{k-1}<\ell\leq r_k.
 \end{cases}
\end{equation}
Since $r_k\geq1$, the stage $\eta_1=y_f$ is always retained.
The reset also satisfies
\begin{equation}\label{eq:resize-contraction}
 \|\eta(t_{k+1}^+)\|^2\leq\|\eta(t_{k+1}^-)\|^2.
\end{equation}

An update is feasible if
\begin{equation}\label{eq:estimated-node-rank}
 \rank\mathscr P_K(\lambda)=p,\qquad
 \lambda\in\widehat\Lambda_k;
\end{equation}
this condition holds automatically when $\widehat\Lambda_k=\varnothing$.
For a feasible update with $\widehat\Lambda_k\neq\varnothing$,
compute $H_k=[H_{1,k}\ \cdots\ H_{r_k,k}]$ from
\begin{equation}\label{eq:estimated-gain-interpolation}
 \sum_{\ell=1}^{r_k}H_{\ell,k}\kappa^{\ell-1}h_f(\lambda)^\ell
 =\mathscr M(\lambda),\qquad\lambda\in\widehat\Lambda_k.
\end{equation}
Here $\mathscr M$ is given by \eqref{eq:interpolation-target}.
The same Vandermonde argument as in
\eqref{eq:interpolation-Vandermonde} gives a unique real $H_k$.
If $\widehat\Lambda_k=\varnothing$, set $H_k=0_{m\times p}$,
as in \eqref{eq:empty-interpolation-gain}.

To bound gains computed from inaccurate frequency estimates,
choose $H_{\max}>0$ and define
\begin{equation}\label{eq:gain-projection}
 \mathcal S_{H_{\max}}(H)=
 \begin{cases}
 H,&\|H\|_F\leq H_{\max},\\
 H_{\max}H/\|H\|_F,&\|H\|_F>H_{\max}.
 \end{cases}
\end{equation}
Apply the bounded gain
\begin{equation}\label{eq:controller-update}
 \widehat H_k=
 \begin{cases}
 \mathcal S_{H_{\max}}(H_k),&\text{if \eqref{eq:estimated-node-rank} holds},\\
 0_{m\times pr_k},&\text{otherwise}.
 \end{cases}
\end{equation}
The applied gain is
\begin{equation}\label{eq:scheduled-gain}
 \widehat H(t)=
 \begin{cases}
 0_{m\times p},&0\leq t<t_1,\\
 \widehat H_k,&t_{k+1}\leq t<t_{k+2},\quad k\geq0.
 \end{cases}
\end{equation}
Thus $\|\widehat H(t)\|_F\leq H_{\max}$, and every update uses
only observations already collected.

The \textit{adaptive controller} consists of the observer and
filter given by the first two equations of
\eqref{eq:complete-controller}, with $r=r(t)$, and the control law
\begin{equation}\label{eq:frequency-controller}
 u(t)=K\xi(t)+\widehat H(t)\eta(t).
\end{equation}
The order $r(t)$, filter state $\eta$, and gain $\widehat H(t)$
are initialized and updated according to
\eqref{eq:online-filter-order}--\eqref{eq:scheduled-gain},
with $\eta_1(0)=y_f(0)$.

Interconnecting this controller with the equivalent plant
\eqref{eq:equivalent-plant} gives, between successive updates,
\begin{equation}\label{eq:adaptive-closed-loop}
\begin{aligned}
 \dd x_{\mathrm{cl}}&=(\mathsf A_{r,\widehat H}x_{\mathrm{cl}}+\mathsf D_rd)\,\dd t
                          +\mathsf B_{m,r}V_m\,\dd W_m
                          +\mathsf B_{p,r}V_p\,\dd W_p,\\
 e_0&=\mathsf C_{r,\widehat H}x_{\mathrm{cl}}+d,
 \qquad x_{\mathrm{cl}}=\col(\chi,\eta,\bar z),
\end{aligned}
\end{equation}
where $r=r(t)$, $\widehat H=\widehat H(t)$, and the matrices
are defined in \eqref{eq:closed-loop-matrices}.
At each update, $\chi$ and $\bar z$ remain continuous,
while $\eta$ is reset according to \eqref{eq:filter-resizing}.

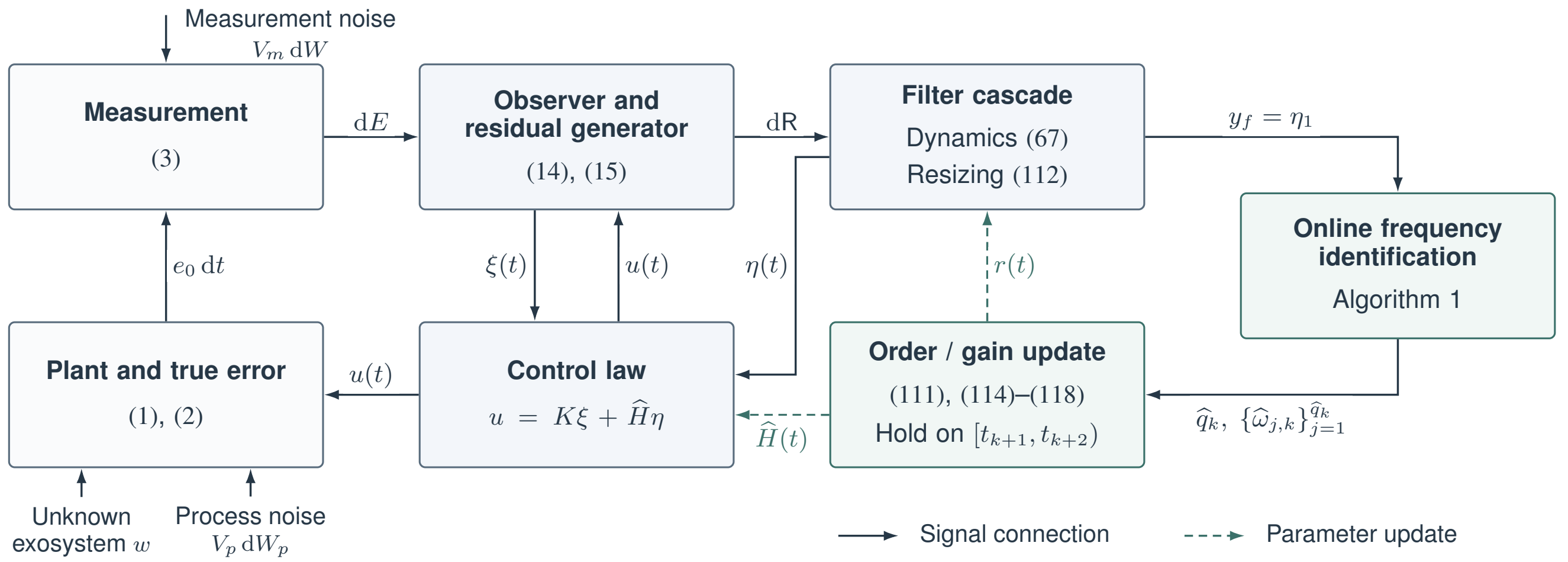
\begin{figure*}[t]
   \centering
   \resizebox{\textwidth}{!}{\begingroup
\definecolor{cflowInk}{HTML}{263746}
\definecolor{cflowBorder}{HTML}{596C7D}
\definecolor{cflowBlue}{HTML}{F3F6FA}
\definecolor{cflowUpdate}{HTML}{3B716A}
\definecolor{cflowGreen}{HTML}{F1F7F5}
\definecolor{cflowPhysical}{HTML}{FAFBFC}
\begin{tikzpicture}[
  font=\sffamily\fontsize{9}{11}\selectfont,
  block/.style={draw=cflowBorder,line width=.65pt,rounded corners=2pt,
    fill=white,text=cflowInk,align=center,text width=3.4cm,
    minimum width=3.9cm,minimum height=1.8cm,inner sep=5pt},
  signal/.style={draw=cflowInk,line width=.72pt,
    -{Latex[length=1.85mm,width=1.25mm]}},
  parameter/.style={draw=cflowUpdate,line width=.72pt,
    dash pattern=on 3pt off 2pt,
    -{Latex[length=1.85mm,width=1.25mm]}},
  label/.style={fill=white,inner sep=2pt,text=cflowInk},
  plabel/.style={fill=white,inner sep=2pt,text=cflowUpdate},
  external/.style={align=center,text=cflowInk,
    font=\sffamily\fontsize{8.5}{10.5}\selectfont}
]

\node[block,fill=cflowPhysical] (measurement) at (0,1.6) {
  \textbf{Measurement}\\[5pt]
  \eqref{eq:measurement-process}
};
\node[block,fill=cflowPhysical] (plant) at (0,-1.6) {
  \textbf{Plant and true error}\\[5pt]
  \eqref{eq:original-system-feedthrough}, \eqref{eq:true-error}
};
\node[block,fill=cflowBlue] (reconstruction) at (5.1,1.6) {
  \textbf{Observer and}\\[-1pt]
  \textbf{residual generator}\\[4pt]
  \eqref{eq:implementable-observer}, \eqref{eq:integrated-residual}
};
\node[block,fill=cflowBlue] (control) at (5.1,-1.6) {
  \textbf{Control law}\\[5pt]
  $u=K\xi+\widehat H\eta$
};
\node[block,fill=cflowBlue] (filter) at (10.2,1.6) {
  \textbf{Filter cascade}\\[4pt]
  Dynamics \eqref{eq:filter-realization}\\[2pt]
  Resizing \eqref{eq:filter-resizing}
};
\node[block,fill=cflowGreen,draw=cflowUpdate] (update) at (10.2,-1.6) {
  \textbf{Order / gain update}\\[4pt]
  \eqref{eq:online-filter-order},
  \eqref{eq:estimated-node-rank}--\eqref{eq:scheduled-gain}\\[3pt]
  Hold on $[t_{k+1},t_{k+2})$
};
\node[block,fill=cflowGreen,draw=cflowUpdate] (identification) at (15.3,0) {
  \textbf{Online frequency}\\[-1pt]
  \textbf{identification}\\[4pt]
  Algorithm~\ref{alg:frequency-identification}
};

\coordinate (noise) at ($(measurement.north)+(0,.60)$);
\node[external] (exo) at (-1.05,-3.30) {
  Unknown\\[-1pt]exosystem $w$
};
\node[external] (processnoise) at (1.05,-3.30) {
  Process noise\\[-1pt]$V_p\,\mathrm dW_p$
};
\draw[signal] (noise) --
  node[external,right=3pt,pos=.40]
    {Measurement noise\\ $V_m\,\mathrm dW$}
  (measurement.north);
\draw[signal] (exo.north) -- ($(plant.south)+(-1.05,0)$);
\draw[signal] (processnoise.north) -- ($(plant.south)+(1.05,0)$);

\draw[signal] (plant.north) --
  node[label,right] {$e_0\,\mathrm dt$} (measurement.south);
\draw[signal] (measurement.east) --
  node[label,above] {$\mathrm dE$} (reconstruction.west);
\draw[signal] (control.west) --
  node[label,above] {$u(t)$} (plant.east);

\draw[signal] ($(reconstruction.south)+(-.52,0)$) --
  node[label,left] {$\xi(t)$} ($(control.north)+(-.52,0)$);
\draw[signal] ($(control.north)+(.52,0)$) --
  node[label,right] {$u(t)$} ($(reconstruction.south)+(.52,0)$);

\draw[signal] (reconstruction.east) --
  node[label,above] {$\mathrm d\mathsf R$} (filter.west);
\draw[signal] ($(filter.west)+(0,-.25)$) -- (7.82,1.35)
  -- (7.82,-1.35) node[label,left,pos=.50] {$\eta(t)$}
  -- ($(control.east)+(0,.25)$);

\draw[signal] (filter.east) -- (15.3,1.6)
  node[label,above,pos=.50] {$y_f=\eta_1$}
  -- (identification.north);
\draw[signal] (identification.south) -- (15.3,-1.6)
  -- (update.east)
  node[label,below,pos=.50]
    {$\widehat q_k,\ \{\widehat\omega_{j,k}\}_{j=1}^{\widehat q_k}$};

\draw[parameter] (update.north) --
  node[plabel,right] {$r(t)$} (filter.south);
\draw[parameter] ($(update.west)+(0,-.25)$) --
  node[plabel,below] {$\widehat H(t)$} ($(control.east)+(0,-.25)$);

\begin{scope}[xshift=5.3cm]
  \draw[signal] (3.05,-3.35) -- (3.80,-3.35);
  \node[anchor=west,text=cflowInk,
        font=\sffamily\fontsize{8.5}{10}\selectfont]
    at (3.94,-3.35) {Signal connection};
\end{scope}
\begin{scope}[xshift=3.5cm]
\draw[parameter] (9.15,-3.35) -- (9.90,-3.35);
\node[anchor=west,text=cflowInk,font=\sffamily\fontsize{8.5}{10}\selectfont]
  at (10.04,-3.35) {Parameter update}; \end{scope}
\end{tikzpicture}
\endgroup}
   \caption{Adaptive error-feedback regulation with online frequency
   identification. Solid arrows denote signal connections; dashed arrows
   denote updates of the filter order and control gain. }
   \label{fig:adaptive-regulation-architecture}
 \end{figure*}

The residual dynamics \eqref{eq:residual-system} remain
independent of the feedback input.
Moreover, the first filter equation in
\eqref{eq:controller-filter} coincides with
\eqref{eq:lowpass-implementation}, and the resizing rule
\eqref{eq:filter-resizing} leaves $\eta_1$ unchanged.
Together with $\eta_1(0)=y_f(0)$, these properties imply
$\eta_1(t)=y_f(t)$ for all $t\geq0$.
Thus the observation model \eqref{eq:identification-model}
remains valid throughout adaptation. The feedback diagram is shown in Fig.~\ref{fig:adaptive-regulation-architecture}.

\subsection{Uniform bounds during order changes}
The fixed-matrix stability in
\eqref{eq:fixed-closed-loop-spectrum} does not by itself bound
the changing-order system \eqref{eq:adaptive-closed-loop}.
The following estimate is uniform over all finite selected orders.
Fourth moments are included to handle the correlation between
$\widehat H$ and $\eta$ in the subsequent comparison.

\begin{proposition}[Bounds under changes of order]
\label{prop:closed-loop-boundedness}
Consider \eqref{eq:adaptive-closed-loop} with Hurwitz $A_o,A_K$,
$\alpha>0$, $0<\kappa<1$, constant $V_m,V_p$, and bounded
deterministic $d$.
At deterministic update times with no finite accumulation,
select a finite order $r\geq1$ and a gain satisfying
$\|\widehat H\|_F\leq H_{\max}<\infty$, both measurable with
respect to the observations available then.
Hold them constant between updates and apply
\eqref{eq:filter-resizing} at every update.

For finite deterministic initial states and order, the closed-loop system \eqref{eq:adaptive-closed-loop}, with
the resizing rule \eqref{eq:filter-resizing}, has a unique
global solution. Moreover, for $a=2,4$,
\begin{equation}\label{eq:closed-loop-moment-bound}
 \sup_{t\geq0}\E\bigl[
 \|x_{\mathrm{cl}}(t)\|^a+\|\xi(t)\|^a+\|u(t)\|^a+\|e_0(t)\|^a
 \bigr]<\infty.
\end{equation}
The bounds are independent of subsequent choices of $r(t)$ in \eqref{eq:online-filter-order}.
\end{proposition}
\begin{IEEEproof}
Between updates, \eqref{eq:adaptive-closed-loop} is a linear SDE
with a finite, fixed order and coefficients determined at the
preceding update. Its unique strong solution, followed by the
reset \eqref{eq:filter-resizing}, defines the solution recursively.
Since the update times have no finite accumulation, this
construction gives a unique global solution.

\emph{A Lyapunov function uniform in $r$:}
Since $A_o$ and $A_K$ are Hurwitz, let $P_o,P_K\succ0$ solve
\begin{equation}\label{eq:fixed-Lyapunov-matrices}
 A_o^\top P_o+P_oA_o=-I_n,\qquad
 A_K^\top P_K+P_KA_K=-I_n.
\end{equation}
With $b=\alpha(1-\kappa)>0$ from
\eqref{eq:uniform-filter-dissipation}, choose $a_\eta>0$
and then $a_\chi>0$ so that
\begin{equation}\label{eq:Lyapunov-weight-choice}
\begin{aligned}
 \delta_\eta&:=a_\eta b-4\|P_KB\|^2H_{\max}^2>0,\\
 \delta_\chi&:=a_\chi-4\|P_KBK\|^2
                     -2a_\eta\alpha^2\|C\|^2/b>0.
\end{aligned}
\end{equation}
For each finite order $r\geq1$ in
\eqref{eq:adaptive-closed-loop}, including those selected by
\eqref{eq:online-filter-order}, define
\begin{equation}\label{eq:variable-order-Lyapunov}
 V_r(x_{\mathrm{cl}})=a_\chi\chi^\top P_o\chi+a_\eta\|\eta\|^2
                                      +\bar z^\top P_K\bar z.
\end{equation}
Let
\[
\begin{aligned}
 c_-&:=\min\{a_\chi\lambda_{\min}(P_o),a_\eta,
                                      \lambda_{\min}(P_K)\},\\
 c_+&:=\max\{a_\chi\lambda_{\max}(P_o),a_\eta,
                                      \lambda_{\max}(P_K)\}
\end{aligned}
\]
 which are independent of $r$ and satisfy
$c_-\|x_{\mathrm{cl}}\|^2\leq V_r(x_{\mathrm{cl}})\leq c_+\|x_{\mathrm{cl}}\|^2$.
Set $V(t):=V_{r(t)}(x_{\mathrm{cl}}(t))$. At every update time $s$,
$\chi$ and $\bar z$ are continuous, whereas
\eqref{eq:resize-contraction} gives
\begin{equation}\label{eq:Lyapunov-jumps}
 V(s^+)-V(s^-)
 =a_\eta\bigl(\|\eta(s^+)\|^2-\|\eta(s^-)\|^2\bigr)
 \leq0.
\end{equation}

\emph{Drift bounds between updates:}
Let $\mathcal L$ denote the It\^o drift operator for
\eqref{eq:adaptive-closed-loop} between updates, and put
$d_\infty:=\sup_{t\geq0}\|d(t)\|<\infty$.
The It\^o differential of \eqref{eq:variable-order-Lyapunov}
along \eqref{eq:adaptive-closed-loop} is
\[
\begin{aligned}
 \dd V_r&=\mathcal LV_r\,\dd t+g_V\,\dd\mathsf W,\\
 \mathsf W&=\col(W_m,W_p),\qquad g_V=[g_m\ g_p],\\
 g_m&=2(-a_\chi\chi^\top P_oL_o
            +a_\eta\alpha\eta_1^\top)V_m,\\
 g_p&=2(a_\chi\chi^\top P_o+\bar z^\top P_K)V_p.
\end{aligned}
\]
Using \eqref{eq:fixed-Lyapunov-matrices} and
\eqref{eq:variable-filter-matrices}, the drift coefficient is
\begin{equation}\label{eq:Lyapunov-drift-expansion}
\begin{aligned}
 \mathcal LV_r={}&
 -a_\chi\|\chi\|^2-2a_\eta\alpha\|\eta\|^2-\|\bar z\|^2+2a_\eta\alpha\kappa
       \sum_{\ell=2}^{r}\eta_\ell^\top\eta_{\ell-1}\\
 &+2a_\eta\alpha\eta_1^\top C\chi
       +2a_\eta\alpha\eta_1^\top d-2\bar z^\top P_KBK\chi \\
 &
       +2\bar z^\top P_KB\widehat H\eta+2\bar z^\top P_KL_od+c_\Sigma,
\end{aligned}
\end{equation}
where, with $\Sigma_m=V_mV_m^\top$ and $\Sigma_p=V_pV_p^\top$,
\begin{equation}\label{eq:Lyapunov-noise-contribution}
\begin{aligned}
 c_\Sigma
 :={}&a_\chi\operatorname{tr}(P_oL_o\Sigma_mL_o^\top)
   +a_\eta\alpha^2\operatorname{tr}(\Sigma_m)\\
 &+\operatorname{tr}\bigl[(a_\chi P_o+P_K)\Sigma_p\bigr].
\end{aligned}
\end{equation}

The estimate \eqref{eq:uniform-filter-dissipation} yields
\[
\begin{aligned}
 &-2a_\eta\alpha\|\eta\|^2
   +2a_\eta\alpha\kappa
      \sum_{\ell=2}^{r}\eta_\ell^\top\eta_{\ell-1}\\
 &\qquad
 =a_\eta\eta^\top
    (\mathsf A_{f,r}+\mathsf A_{f,r}^\top)\eta
 \leq-2a_\eta b\|\eta\|^2.
\end{aligned}
\]
By Young's inequality, $\|\eta_1\|\leq\|\eta\|$, and
$\|\widehat H\|\leq\|\widehat H\|_F\leq H_{\max}$,
\begin{equation}\label{eq:filter-cross-term-bounds}
\begin{aligned}
 2a_\eta\alpha|\eta_1^\top C\chi|
 &\leq\tfrac12a_\eta b\|\eta\|^2
       +\frac{2a_\eta\alpha^2\|C\|^2}{b}\|\chi\|^2,\\
 2a_\eta\alpha|\eta_1^\top d|
 &\leq\tfrac12a_\eta b\|\eta\|^2
       +\frac{2a_\eta\alpha^2}{b}d_\infty^2,
\end{aligned}
\end{equation}
and
\begin{equation}\label{eq:plant-cross-term-bounds}
\begin{aligned}
 2|\bar z^\top P_KBK\chi|
 &\leq\tfrac14\|\bar z\|^2
       +4\|P_KBK\|^2\|\chi\|^2,\\
 2|\bar z^\top P_KB\widehat H\eta|
 &\leq\tfrac14\|\bar z\|^2
       +4\|P_KB\|^2H_{\max}^2\|\eta\|^2,\\
 2|\bar z^\top P_KL_od|
 &\leq\tfrac14\|\bar z\|^2
       +4\|P_KL_o\|^2d_\infty^2.
\end{aligned}
\end{equation}
Applying these bounds to \eqref{eq:Lyapunov-drift-expansion},
with the choices in \eqref{eq:Lyapunov-weight-choice}, yields
\begin{equation}\label{eq:uniform-Lyapunov-drift}
\begin{aligned}
 \mathcal LV_r
 &\leq-\delta_\chi\|\chi\|^2-\delta_\eta\|\eta\|^2
       -\tfrac14\|\bar z\|^2+C_0\\
 &\leq-c_0V_r+C_0,
\end{aligned}
\end{equation}
where $c_0:=\frac{\min\{\delta_\chi,\delta_\eta,1/4\}}{c_+}>0, $ and
\[
\begin{aligned}
 C_0:=\left(\frac{2a_\eta\alpha^2}{b}
              +4\|P_KL_o\|^2\right)d_\infty^2+c_\Sigma.
\end{aligned}
\]
These constants do not depend on $r$ or the update times.

The bound $c_-\|x_{\mathrm{cl}}\|^2\leq V_r$ implies
\[
\begin{aligned}
 \|g_V\|^2&=\|g_m\|^2+\|g_p\|^2\\
 &\leq8\|V_m\|^2
   \bigl(a_\chi^2\|P_oL_o\|^2\|\chi\|^2
         +a_\eta^2\alpha^2\|\eta\|^2\bigr)\\
 &\quad+8\|V_p\|^2
   \bigl(a_\chi^2\|P_o\|^2\|\chi\|^2
         +\|P_K\|^2\|\bar z\|^2\bigr)\\
 &\leq C_gV_r
\end{aligned}
\]
for some $C_g\geq0$ independent of $r$. Applying It\^o's
formula to $V_r^2$ and using \eqref{eq:uniform-Lyapunov-drift}
therefore gives
\begin{equation}\label{eq:fourth-moment-drift}
\begin{aligned}
 \mathcal L(V_r^2)
 &=2V_r\mathcal LV_r+\|g_V\|^2\leq-2c_0V_r^2+(2C_0+C_g)V_r\\
 &\qquad \leq-c_0V_r^2+C_2,
\end{aligned}
\end{equation}
where $C_2:=(2C_0+C_g)^2/(4c_0)$.

\emph{Moment bounds across updates:}
For $N>V(0)$, set
$\tau_N:=\inf\{t\geq0:V(t)\geq N\}$.
Since $V(s)<N$ for $s<\tau_N$ and $\|g_V\|^2\leq C_gV$,
the processes
\[
\begin{gathered}
 \int_0^{t\wedge\tau_N}e^{c_0s}g_V(s)\,\dd\mathsf W(s),\ 
 \int_0^{t\wedge\tau_N}2e^{c_0s}V(s)g_V(s)\,\dd\mathsf W(s)
\end{gathered}
\]
are square-integrable martingales on finite intervals
and have zero mean.
Applying It\^o's formula to $e^{c_0t}V(t)$ and
$e^{c_0t}V(t)^2$ up to $t\wedge\tau_N$, taking expectations,
and using \eqref{eq:uniform-Lyapunov-drift},
\eqref{eq:fourth-moment-drift}, and the nonpositive jumps
in \eqref{eq:Lyapunov-jumps}, we obtain
\[
\begin{aligned}
 \E\!\left[e^{c_0(t\wedge\tau_N)}V(t\wedge\tau_N)\right]
 &\leq V(0)+\frac{C_0}{c_0}(e^{c_0t}-1),\\
 \E\!\left[e^{c_0(t\wedge\tau_N)}V(t\wedge\tau_N)^2\right]
 &\leq V(0)^2+\frac{C_2}{c_0}(e^{c_0t}-1).
\end{aligned}
\]
Since $\tau_N\to\infty$ a.s., Fatou's lemma \cite[Thm.~1.6.5]{Durrett2019} yields
\begin{equation}\label{eq:residual-fourth-moments}
\begin{aligned}
 \E V(t)&\leq e^{-c_0t}V(0)
                       +\frac{C_0}{c_0}(1-e^{-c_0t}),\\
 \E[V(t)^2]&\leq e^{-c_0t}V(0)^2
                       +\frac{C_2}{c_0}(1-e^{-c_0t}).
\end{aligned}
\end{equation}
Together with $c_-\|x_{\mathrm{cl}}\|^2\leq V$, these estimates bound
$\E\|x_{\mathrm{cl}}(t)\|^2$ and $\E\|x_{\mathrm{cl}}(t)\|^4$ uniformly in $t$.
Moreover, $\xi=\bar z-\chi$, \eqref{eq:frequency-controller},
and the output equation in \eqref{eq:fixed-closed-loop} give
\[
 \|\xi\|+\|u\|+\|e_0\|
 \leq c_{\mathrm{out}}\|x_{\mathrm{cl}}\|+d_\infty
\]
for a constant $c_{\mathrm{out}}$ independent of $r$.
This proves \eqref{eq:closed-loop-moment-bound}.
\end{IEEEproof}
\subsection{Adaptive regulation and limiting mean-square error}

\begin{assumption}[Gain budget]\label{ass:gain-budget}
The threshold $H_{\max}$ satisfies
\begin{equation}\label{eq:gain-budget}
 \|H_\star\|_F<H_{\max},
\end{equation}
where $H_\star$ is defined by \eqref{eq:ideal-gain}.
\end{assumption}

\begin{remark}
	For each fixed admissible set of active frequencies, the ideal
gain $H_\star$ defined by \eqref{eq:ideal-gain} has finite norm.
Hence Assumption~\ref{ass:gain-budget} can always be satisfied
by choosing $H_{\max}$ sufficiently large. During identification, inaccurate estimates may lead to
controller gains that are too large and strongly amplify noise.
The gain limit keeps these gains within a fixed bound and
ensures the moment bounds in
Proposition~\ref{prop:closed-loop-boundedness}.
Assumption~\ref{ass:gain-budget} ensures that this bound is
still large enough to include the ideal gain $H_\star$.
As the estimates converge, the projection
\eqref{eq:gain-projection} eventually stops modifying the
computed gain, allowing the adaptive controller to attain
the same asymptotic performance as the known-frequency controller.
\end{remark}

For comparison, use the fixed controller $(r_\star,H_\star)$
in \eqref{eq:ideal-gain}, driven by the same $d,W_m,W_p$ and with
the same initial plant and observer states as the adaptive system.
Initialize its filter by
$\eta_{\star,1}(0)=y_f(0)$ and
$\eta_{\star,\ell}(0)=0$ for $\ell\geq2$.
When comparing gains or filter states of different dimensions,
append zero blocks to the shorter object.

\begin{theorem}[Adaptive output regulation]
\label{thm:adaptive-regulation}
Suppose that Assumptions~\ref{ass:plant-exosystem}--
\ref{ass:gain-budget} hold, and choose
$K,L_o,\alpha,\kappa$ as in
Theorem~\ref{thm:known-frequency-regulation}.
Then Algorithm~\ref{alg:frequency-identification},
together with the adaptive controller
\eqref{eq:frequency-controller}, solves
Problem~\ref{prob:adaptive-regulation}.

Moreover, for the known-frequency comparison system
\eqref{eq:fixed-closed-loop} with $(r_\star,H_\star)$
in \eqref{eq:ideal-gain},
\begin{equation}\label{eq:ideal-comparison-convergence}
\begin{aligned}
 \E\|\bar z(t)-\bar z_\star(t)\|^2\rightarrow0,\quad
 \E\|e_0(t)-e_{0,\star}(t)\|^2\rightarrow0,
\end{aligned}
\end{equation}
 the true tracking error satisfies
\begin{equation}\label{eq:stochastic-regulation}
\begin{aligned}
 \E[e_0(t)]\rightarrow0,\quad
 \E[e_0(t)e_0(t)^\top]\rightarrow \Sigma_e,
\end{aligned}
\end{equation}
as $t\to\infty$, and
\begin{equation}\label{eq:adaptive-mse-limit}
\begin{aligned}
 \lim_{t\to\infty}\E\|e_0(t)\|^2
 &=J_m+J_p=J\\
 \leq &\|\mathcal M_{e,m}\|\operatorname{tr}(\Sigma_m)
       +\|\mathcal M_{e,p}\|\operatorname{tr}(\Sigma_p),
\end{aligned}
\end{equation}
where $\Sigma_e,J$ and $J_m,J_p$ are defined by
\eqref{eq:limiting-error-covariance} and \eqref{eq:noise-performance-functional}.
In the absence of both noises, $e_0(t)\to0$.
\end{theorem}
\begin{IEEEproof}
Since $d$ is bounded and
$\|\widehat H(t)\|_F\leq H_{\max}$ by
\eqref{eq:controller-update}--\eqref{eq:scheduled-gain},
Proposition~\ref{prop:closed-loop-boundedness} gives global
well-posedness and \eqref{eq:closed-loop-moment-bound}.
The residual equation \eqref{eq:residual-system} is independent
of $u$, and the resizing rule \eqref{eq:filter-resizing} leaves
$\eta_1$ unchanged. With the prescribed initializations,
uniqueness in \eqref{eq:residual-system} and
\eqref{eq:controller-filter} therefore gives
\begin{equation}\label{eq:shared-residual-states}
 \chi_\star=\chi,\qquad \eta_{\star,1}=\eta_1=y_f.
\end{equation}
Thus the observation model \eqref{eq:identification-model}
remains valid under feedback.

\emph{Order and gain convergence:}
By Theorem~\ref{thm:identification} and
\eqref{eq:online-filter-order}, almost surely,
$r_k=r_\star$ eventually and the estimated nodes converge
to the true nodes. The rank condition
\eqref{eq:known-node-rank} persists at nearby nodes, and
\eqref{eq:interpolation-target} and
\eqref{eq:interpolation-Vandermonde} show that the interpolation
gain depends continuously on these nodes. Hence
\eqref{eq:estimated-gain-interpolation} gives $H_k\to H_\star$.

By \eqref{eq:gain-budget} and $H_k\to H_\star$,
$\|H_k\|_F<H_{\max}$ for all sufficiently large $k$.
Since \eqref{eq:estimated-node-rank} also holds eventually,
\eqref{eq:gain-projection} and \eqref{eq:controller-update}
give $\widehat H_k=H_k$ for all sufficiently large $k$.
The holding rule \eqref{eq:scheduled-gain} therefore yields
$\widehat H(t)\to H_\star$ a.s. For every $1\leq a<\infty$,
\[
 \|\widehat H(t)-H_\star\|_F^a
 \leq (H_{\max}+\|H_\star\|_F)^a.
\]
Together with the almost-sure convergence proved above, the dominated convergence theorem  implies
\begin{equation}\label{eq:gain-consistency}
\begin{aligned}
 \E\|\widehat H(t)-H_\star\|_F^a&\longrightarrow0,
 \qquad 1\leq a<\infty.
\end{aligned}
\end{equation}

\emph{Filter and correction-input convergence:}
Almost surely, there is a finite update time $\tau$ such that
$r(t)=r_\star$ for all $t\geq\tau$.
By \eqref{eq:filter-resizing}, $\eta$ has no jumps at updates
after $\tau$. Write $\tilde{\eta}=\eta-\eta_\star$, appending zero
components to the shorter vector when the dimensions differ.
Subtracting the filter equations in
\eqref{eq:fixed-closed-loop} and using
\eqref{eq:shared-residual-states}, we obtain
\(
 \dot{\tilde{\eta}}=\mathsf A_{f,r_\star}\tilde{\eta},\) for \(
 t>\tau,
\)
and so
\begin{equation}\label{eq:filter-comparison}
 \tilde{\eta}(t)
 =e^{\mathsf A_{f,r_\star}(t-\tau)}
      \tilde{\eta}(\tau^+),\qquad t\geq\tau.
\end{equation}
By \eqref{eq:uniform-filter-dissipation},
$\tilde{\eta}(t)\to0$ a.s. The fourth-moment bounds
\eqref{eq:closed-loop-moment-bound} and
\eqref{eq:known-frequency-moment-bound} imply
$\sup_{t\geq0}\E\|\tilde{\eta}(t)\|^4<\infty$.
Thus, for every $\varepsilon>0$, 
using the Cauchy--Schwarz inequality, we have
\[
\begin{aligned}
 \E\|\tilde{\eta}(t)\|^2
 &\leq \varepsilon^2
   +\E\!\left[
      \|\tilde{\eta}(t)\|^2
      \mathbf{1}_{\{\|\tilde{\eta}(t)\|>\varepsilon\}}
   \right]\\
 &\leq \varepsilon^2
   +\bigl(\E\|\tilde{\eta}(t)\|^4\bigr)^{1/2}
      \Prob\{\|\tilde{\eta}(t)\|>\varepsilon\}^{1/2}\\
 &\leq \varepsilon^2
   +c\,\Prob\{\|\tilde{\eta}(t)\|>\varepsilon\}^{1/2},
\end{aligned}
\]
where $\mathbf{1}_A$ denotes the indicator of an event $A$,
and $c$ is independent of $t$ and $\varepsilon$.
Letting $t\to\infty$ and then $\varepsilon\downarrow0$
proves $\E\|\tilde{\eta}(t)\|^2\to0$.

Define
\begin{equation}\label{eq:adaptive-forcing}
\begin{aligned}
 f&:=\widehat H\eta-H_\star\eta_\star
  =(\widehat H-H_\star)\eta+H_\star\tilde{\eta}.
\end{aligned}
\end{equation}
Using Cauchy--Schwarz inequality, together with \eqref{eq:gain-consistency} and
\eqref{eq:closed-loop-moment-bound}, we obtain that, as  $t\to \infty$,
\begin{equation}\label{eq:adaptive-forcing-L2}
\begin{aligned}
 \E\|f(t)\|^2
 &\leq 2\E\!\left[
       \|\widehat H(t)-H_\star\|_F^2\|\eta(t)\|^2
       \right]+2\|H_\star\|^2\E\|\tilde\eta(t)\|^2\\
 &\leq 2\bigl(\E\|\widehat H(t)-H_\star\|_F^4\bigr)^{1/2}
         \bigl(\E\|\eta(t)\|^4\bigr)^{1/2}\\
 &\quad+2\|H_\star\|^2\E\|\tilde\eta(t)\|^2
 \rightarrow0.
\end{aligned}
\end{equation}

\emph{State and error convergence:}
Set $\tilde{z}=\bar z-\bar z_\star$ and
$\tilde{e}=e_0-e_{0,\star}$.
From \eqref{eq:fixed-closed-loop} and
\eqref{eq:shared-residual-states}, the process-noise increments
cancel because both plants are driven by the same $W_p$. Hence
\begin{equation}\label{eq:comparison-dynamics}
\begin{aligned}
 \dot{\tilde{z}}&=A_K\tilde{z}+Bf,\qquad \tilde{z}(0)=0,\\
 \tilde{e}&=C_K\tilde{z}+Df.
\end{aligned}
\end{equation}
Since the plant states are continuous at updates,
\[
 \tilde{z}(t)=\int_0^t e^{A_K(t-s)}Bf(s)\,\dd s.
\]
With $\|X\|_{L^2}:=\|X\|_{L^2(\mathbb P)}=(\E\|X\|^2)^{1/2}$, exponential stability
of $A_K$ and Minkowski's inequality yield, for $t\geq s_0\geq0$,
\begin{equation}\label{eq:comparison-convolution-bound}
\begin{aligned}
 \|\tilde{z}(t)\|_{L^2}
 &\leq c e^{-\beta_K(t-s_0)}
          \int_0^{s_0}\|f(s)\|_{L^2}\,\dd s\\
 &\quad{}+c\sup_{s\geq s_0}\|f(s)\|_{L^2},
\end{aligned}
\end{equation}
where $c,\beta_K>0$ are independent of $t,s_0$.
Moreover, $\sup_{t\geq0}\E\|f(t)\|^2<\infty$, and
$\E\|f(t)\|^2\to0$ by \eqref{eq:adaptive-forcing-L2}.
Letting $t\to\infty$ and then $s_0\to\infty$ therefore gives
$\|\tilde{z}(t)\|_{L^2}\to0$. Moreover,
\begin{equation}\label{eq:conve-tilde-e}
	 \|\tilde{e}(t)\|_{L^2}
 \leq\|C_K\|\|\tilde{z}(t)\|_{L^2}
       +\|D\|\|f(t)\|_{L^2}\rightarrow0,
\end{equation}
which proves \eqref{eq:ideal-comparison-convergence}.

\emph{Limiting error moments:}
The convergence in \eqref{eq:conve-tilde-e} implies
\[
 \|\E[e_0(t)]-\E[e_{0,\star}(t)]\|
 \leq\|\tilde{e}(t)\|_{L^2}\longrightarrow0.
\]
Using $\tilde e=e_0-e_{0,\star}$, we have
\(
 e_0e_0^\top-e_{0,\star}e_{0,\star}^\top
 =\tilde e\,e_0^\top+e_{0,\star}\tilde e^\top.
\)
Taking expectations and applying Cauchy--Schwarz gives
\begin{equation}\label{eq:error-second-moment-comparison}
\begin{aligned}
 &\|\E[e_0e_0^\top]-\E[e_{0,\star}e_{0,\star}^\top]\|\\
 &\quad\leq\|\tilde e\|_{L^2}
       \bigl(\|e_0\|_{L^2}+\|e_{0,\star}\|_{L^2}\bigr)
       \longrightarrow0,
\end{aligned}
\end{equation}
as $t\to \infty$,
since $\|\tilde e(t)\|_{L^2}\to0$,  $\|e_0\|_{L^2}$ and $\|e_{0,\star}\|_{L^2}$
 are uniformly bounded in $t$.
Together with \eqref{eq:known-frequency-regulation}, these
limits prove \eqref{eq:stochastic-regulation}.
Since $\E\|e_0(t)\|^2
=\operatorname{tr}(\E[e_0(t)e_0(t)^\top])$, using \eqref{eq:error-second-moment-comparison} and \eqref{eq:noise-performance-functional}, we obtain
\[
\begin{aligned}
 \lim_{t\to\infty}\E\|e_0(t)\|^2
 &=\operatorname{tr}(\Sigma_e)=J_m+J_p\\
 &\leq\|\mathcal M_{e,m}\|\operatorname{tr}(\Sigma_m)
       +\|\mathcal M_{e,p}\|\operatorname{tr}(\Sigma_p).
\end{aligned}
\]
This proves \eqref{eq:adaptive-mse-limit}.
If $V_m=V_p=0$, the error is deterministic and
\eqref{eq:adaptive-mse-limit} yields $e_0(t)\to0$.
\end{IEEEproof}

\section{Numerical Simulations}
\label{sec:numerical-simulations}
We consider an unstable plant to examine
frequency identification, closed-loop behavior, and the effect of
measurement noise on tracking performance. All simulations below
use $V_p=0$, so $J_p=0$ and $J=J_m$; they retain the
measurement-noise baseline of the proposed method.

\subsection{Simulation setup}
For \eqref{eq:original-system-feedthrough}, choose
\[
\begin{gathered}
 A=\begin{bmatrix}0.2&1\\1&0\end{bmatrix},\quad
 B=\begin{bmatrix}0\\1\end{bmatrix},\quad
 C=[1\;0.5],\quad D=0.1,\\
 P=\begin{bmatrix}
 0&0&0&0&0\\
 0.2&0.3&0&0.15&0
 \end{bmatrix},\quad
 Q=[-1\;-1\;0\;-0.6\;0].
\end{gathered}
\]
The exosystem is given by
\[
\begin{gathered}
 S=\operatorname{diag}\bigl(0,S_{0.8},S_{3.2}\bigr),\qquad
 w(0)=\col(1,1,0,1,0),\\
 S_\omega=\begin{bmatrix}0&-\omega\\\omega&0\end{bmatrix}.
\end{gathered}
\]
The output $y=Cz+Du$ is required to track
$y_{\mathrm{ref}}=-Qw=1+\cos(0.8t)+0.6\cos(3.2t)$
in the presence of the disturbance $Pw$.
The controller receives only the noisy error measurement
\begin{equation}\label{eq:sim-noise}
 \dd E=e_0\,\dd t+\delta\,\dd W_m,
 \qquad e_0=y-y_{\mathrm{ref}},\quad \Sigma_m=\delta^2.
\end{equation}

We choose $K=[-3\;-3]$, $L_o=\col(2,2)$, $\alpha=2$,
$\kappa=0.9$, and $H_{\max}=500$.
Both $A_K$ and $A_o$ are Hurwitz, and the plant has no
imaginary-axis transmission zeros.
The corresponding ideal controller has order $r_\star=5$ and
$\|H_\star\|_F\simeq339.275<H_{\max}$.
The initial conditions are $z(0)=\col(0.5,-0.3)$, $\xi(0)=0$,
and zero filter states, with initial adaptive order $r(0)=1$.

Algorithm~\ref{alg:frequency-identification} uses the band
$[\omega_-,\omega_+]=[0.5,3.5]$, $t_{\mathrm{b}}=0.5$,
and $\gamma_{\mathrm{b}}=0.15$.
The simulation runs on $[0,400]$, with updates at
$1,2,4, 8, \ldots,256$.

\subsection{Identification and regulation}
\begin{figure}[t]
 \centering
 \includegraphics[width=\columnwidth]
 {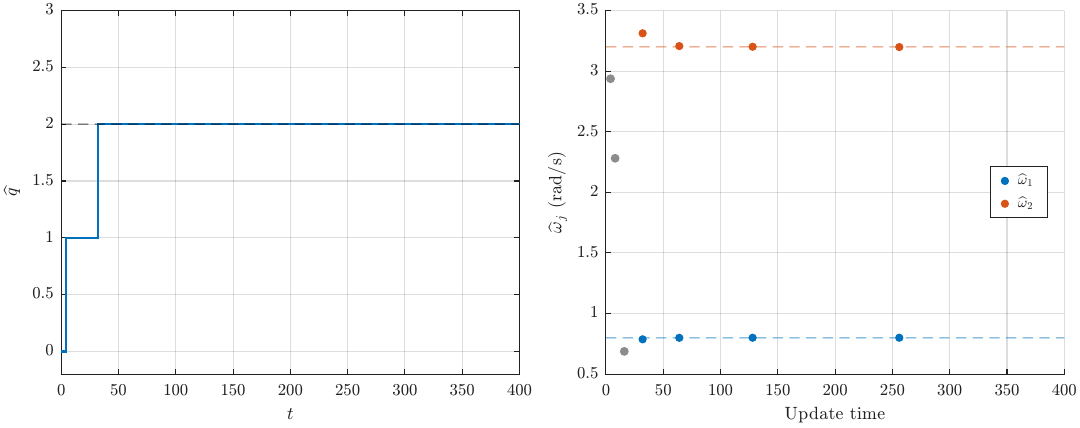}\par
 \caption{Frequency identification for $\delta=0.05$;
 dashed lines mark the true values.}  \label{fig:sim-identification}
 \includegraphics[width=\columnwidth]
 {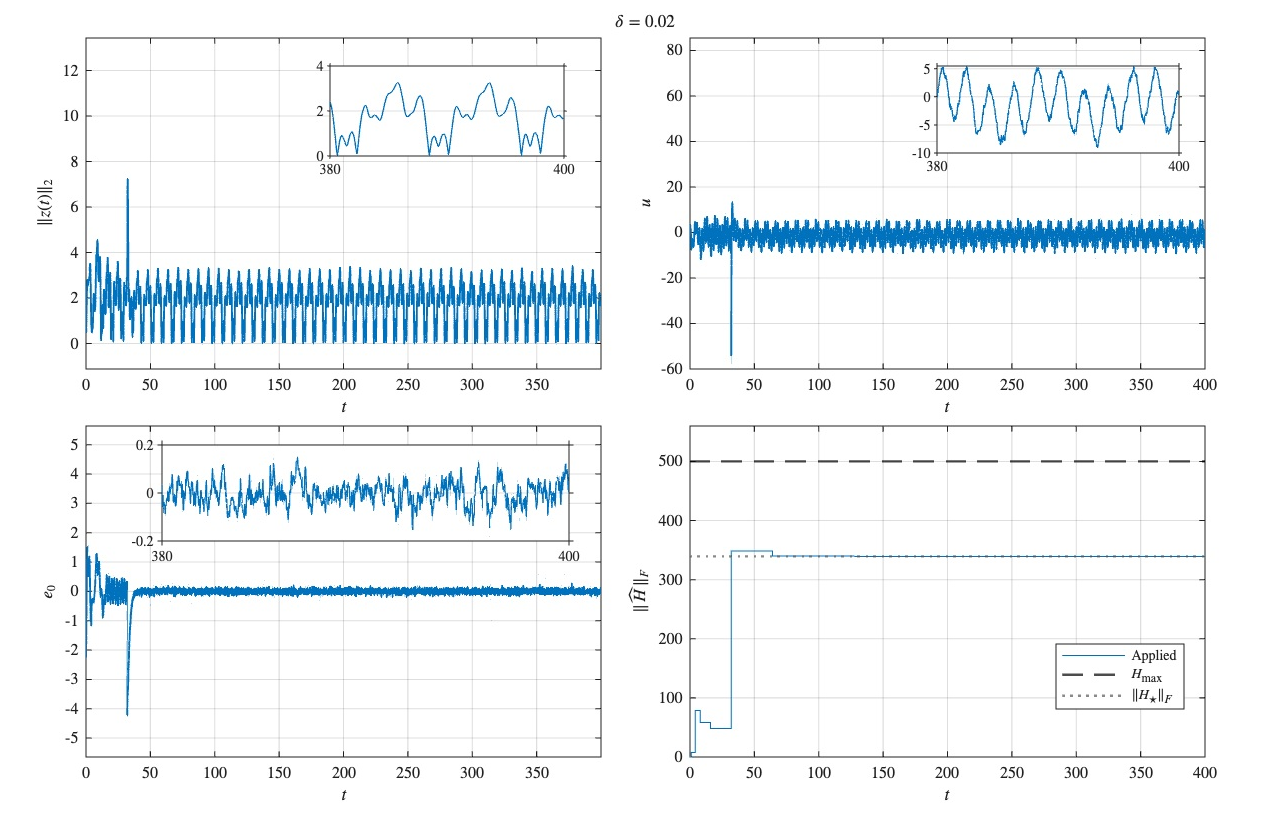}
 \caption{ Closed-loop
 responses for $\delta=0.02$, showing $\|z(t)\|_2$, $u$,
 $e_0$, and $\|\widehat H\|_F$ in row order.
 Insets show late-time details.}
 \label{fig:sim-closed-loop}
\end{figure}

 Fig.~\ref{fig:sim-identification} shows
frequency identification for $\delta=0.05$.
The estimated number of positive frequencies and their values
are compared with $q=2$ and $\omega_1=0.8$, $\omega_2=3.2$.
The estimator uses the filtered signal $y_f$, and its outputs
update the controller order and gain according to
\eqref{eq:online-filter-order}--\eqref{eq:scheduled-gain}.

Fig.~\ref{fig:sim-closed-loop} shows the closed-loop response for $\delta=0.02$,
including the plant state norm $\|z(t)\|_2$, control input $u$,
true tracking error $e_0$, and applied gain norm
$\|\widehat H\|_F$.
The gain panel includes $H_{\max}$ and $\|H_\star\|_F$
as references for the imposed limit and ideal gain, respectively.
The main axes retain the initial transient, while the insets
show late-time oscillations and noise-induced fluctuations.
Persistent oscillations in the state and input are compatible
with regulation because the reference and disturbance remain
time varying.

\subsection{Effect of measurement noise}
To examine how measurement noise affects tracking, we consider
$\delta\in\{0,0.02,0.05,0.2\}$.
For each positive noise level, we repeat the simulation with
independent noise realizations and the same initial conditions.
Each run performs its own frequency estimation and controller
updates. At each time $t$, we average the squared tracking
errors over these runs:
\begin{equation}\label{eq:sim-empirical-mse}
 \widehat J_\delta(t)=\frac{1}{M_\delta}
 \sum_{a=1}^{M_\delta}|e_0^{(a)}(t)|^2,
\end{equation}
where $e_0^{(a)}$ is the tracking error in the $a$th run.
We use $M_\delta=50$ for $\delta>0$ and $M_0=1$, since the
noiseless simulation is deterministic.

For this example, \eqref{eq:noise-performance-functional}
predicts the limiting mean-square error
$
 J_m\simeq7.379465\,\delta^2.
$

\begin{figure}[t]
 \centering
 \includegraphics[width=\columnwidth]
 {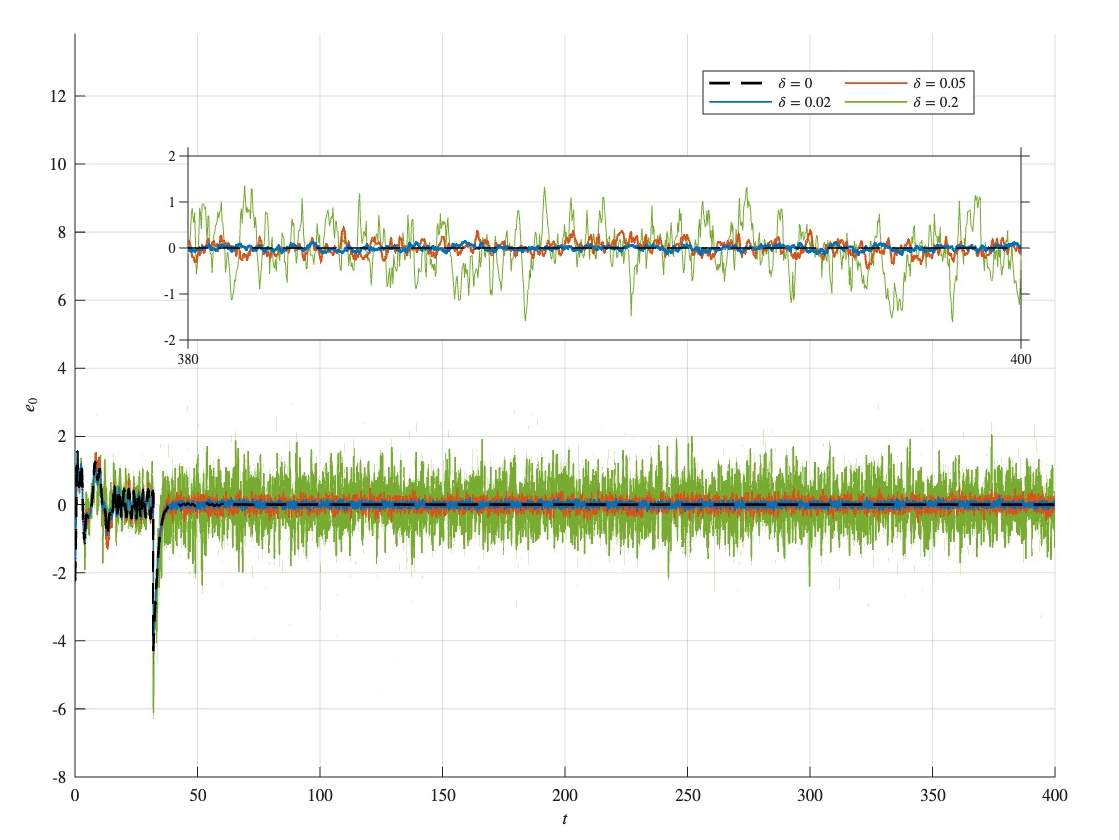}\par
 \includegraphics[width=\columnwidth]
 {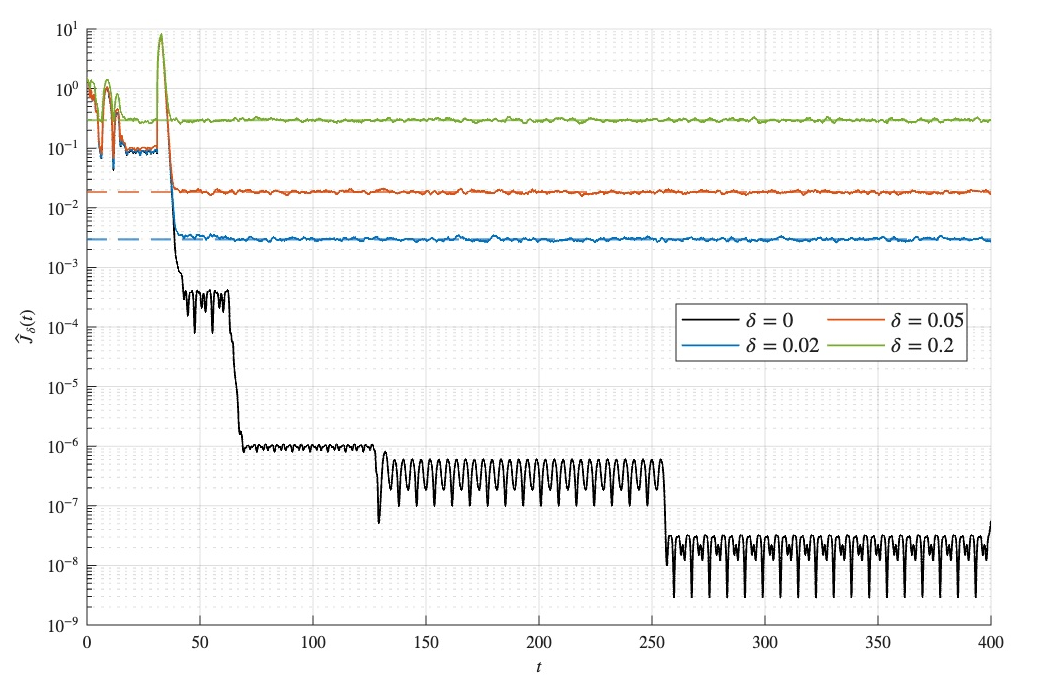}
 \caption{Tracking performance for $\delta=0,0.02,0.05,0.2$.
 Top: representative errors. Bottom: empirical MSE and
 theoretical limits $J_m$ (dashed, $\delta>0$).
}
 \label{fig:sim-noise-traces}
 \label{fig:sim-noise-mse}
\end{figure}

The upper part of Figure~\ref{fig:sim-noise-traces} shows one
tracking-error trajectory for each noise level.
The lower part compares the averaged squared errors
$\widehat J_\delta(t)$ with the theoretical limits $J_m$.

\section{Conclusion}
This paper developed an adaptive error-feedback regulator for
linear systems with an unknown harmonic exosystem and additive
process and measurement noise. An equivalent-system transformation
and an implementable observer provide an input-independent residual
for online frequency identification. The controller uses only the
identified frequency count and values, while uniform moment bounds
are maintained during filter-order changes. Under the stated
assumptions and gain-budget condition, the adaptive plant response
and tracking error approach those of the known-frequency regulator
in mean square. 
Possible future extensions include systems with multiplicative noise and
infinite-dimensional plants. Multiplicative noise calls for further
analysis of the residual statistics and stochastic stability.
For infinite-dimensional systems, particular attention is needed to well-posedness and to possibly unbounded input and observation operators.

\IEEEtriggeratref{33}

\end{document}